\documentclass[12pt,reqno]{amsart}
\usepackage{graphicx}
\usepackage{amsmath, amsthm, amssymb, stmaryrd}
\usepackage{hyperref}
\usepackage{enumerate}
\usepackage{url}
\usepackage{scalerel}

\usepackage{mathtools}
\DeclarePairedDelimiter\floor{\lfloor}{\rfloor}
\DeclarePairedDelimiter\ceil{\lceil}{\rceil}

\usepackage{multirow, array}
\usepackage{placeins}

\usepackage{caption} 
\advance \topmargin by -\headheight
\advance \topmargin by -\headsep
     
\evensidemargin \oddsidemargin
\newtheorem{thm}{Theorem}[section]
\newtheorem{lemma}[thm]{Lemma}
\newtheorem{prop}[thm]{Proposition}
\newtheorem{cor}[thm]{Corollary}

\theoremstyle{definition}

\theoremstyle{remark}

\numberwithin{equation}{section}

\newcommand*\wrapletters[1]{\wr@pletters#1\@nil}
\def\wr@pletters#1#2\@nil{#1\allowbreak\if&#2&\else\wr@pletters#2\@nil\fi}

\def\eps{\varepsilon}

\def\le{\leqslant} \def\ge{\geqslant}

\def \bC {\mathbb C}

\def \bN {\mathbb N}

\def \bQ {\mathbb Q}
\def \bR {\mathbb R}
\def \bZ {\mathbb Z}

\def \bc {\mathbf c}

\def \bh {\mathbf h}

\def \bx {\mathbf x}

\def \fm {\mathfrak m}

\def \fM {\mathfrak M}

\def \lcm {\mathrm{lcm}}

\def \deg {\mathrm{deg}}

\def \dens {{\mathrm{dens}}}

\title[On the probability that the integer parts of polynomials are coprime]{On the probability that the integer parts of polynomials are coprime}

\author[Philippa Holdridge]{Philippa Holdridge}
\address{HUN-REN Alfr\'ed R\'enyi Institute of Mathematics, Re\'altanoda utca 13--15., H-1053 Budapest,  Hungary;   \newline \hspace*{4mm}
MTA--HUN-REN RI Lend\"ulet ``Momentum'' Arithmetic Combinatorics Research Group, Re\'altanoda utca 13--15., H-1053 Budapest,  Hungary}
\email{holdridge.philippa@renyi.hu}

\thanks{}
\date{}

\begin{document}

\begin{abstract}
We prove that for any $s\ge 2$ and polynomials $f_{1}(x),\dots,f_{s-1}(x)\in \bR[x]$ satisfying a certain equidistribution condition, the natural density of the numbers $n\in \bN$ such that $n,\floor{f_{1}(n)},\dots,\floor{f_{s-1}(n)}$ are mutually coprime is $1/\zeta(s)$.

This generalises a recent result of Champagne, in which this was proven under the assumption that $f_{1}$ is linear.
\end{abstract}

\maketitle

\setcounter{tocdepth}{1}
\tableofcontents

\section{Introduction}
It is well-known that if $s\ge 2$ is an integer and $n_{1},\dots,n_{s}$ are random variables chosen uniformly and independently on $\{1,2,\dots,N\}$, then the probability that $\gcd(n_{1},\dots,n_{s})=1$ tends to $1/\zeta(s)$ as $N\rightarrow\infty$. Here $\zeta$ denotes the Riemann zeta function. Instead of choosing these random integers independently, we can ask how often $\gcd(n,\floor{f_{1}(n)},\dots,\floor{f_{s-1}(n)})=1$, where $n$ is chosen uniformly on $\{1,2,\dots,N\}$ and $f_{1},\dots,f_{s-1}:\bN\rightarrow \bR$ are fixed functions. For certain choices of functions $f_{1},\dots,f_{s-1}$, we might expect that the probability still tends to $1/\zeta(s)$.

For $A\subseteq \bN$, define the upper and lower natural densities by
\[\overline{\dens} A=\limsup_{N\rightarrow\infty}\frac{1}{N}\#\{n\in A:n\le N\},\]
\[\underline{\dens} A=\liminf_{N\rightarrow\infty}\frac{1}{N}\#\{n\in A:n\le N\}\]
and whenever $\overline{\dens} A=\underline{\dens} A$, we define the natural density as $\dens A=\overline{\dens} A$. Then we may state our problem as that of finding certain classes of functions $f_{1},\dots,f_{s-1}$ for which
\begin{equation}\label{main_eq}\dens \{n\in \bN: \gcd(n,\floor{f_{1}(n)},\dots,\floor{f_{s-1}(n)})=1\}=\frac{1}{\zeta(s)}.\end{equation}

Of particular interest to us is the case when $f_{1},\dots,f_{s-1}$ are polynomials. If $f_{i}(n)=\alpha_{j}n^{m_{j}}$ for some irrational numbers $\alpha_{1},\dots,\alpha_{s-1}$ and integers $1=m_{1}<\cdots<m_{s-1}$, then Banks and Shparlinski \cite{bs2025} proved \eqref{main_eq} under the assumption that the $\alpha_{1},\dots,\alpha_{s-1}$ are not too well-approximated by rationals. More precisely, if we define the exponential type of a number $\alpha$ to be
\[\sup\{\theta\in \bR: \text{for infinitely many }q\in \bN, \:\|q\alpha\|\le \exp(-q^{\theta})\}.\]
Then we require that all of the $\alpha_{i}$ have finite exponential type. Numbers of infinite exponential type are very rare, for example, such numbers must be Liouville and therefore form a set of Hausdorff dimension $0$. But there are still uncountably many of them. For example, if $a_{n}$ is any strictly increasing sequence of integers such that $a_{n+1}\ge 2^{2^{a_{n}}}$ for all $n$, then
\[\alpha=\sum_{n=1}^{\infty} 2^{-a_{n}}\]
has infinite exponential type. For $q=2^{a_{n}}$, we have $\|q\alpha\|\ll 2^{-a_{n+1}}$.

More recently, Champagne, \cite{c2026} was able to remove the finite exponential type condition from Banks and Shparlinski's result, and more generally was able to prove \eqref{main_eq} for a wider range of $f_{1},\dots,f_{s-1}$. Specifically, it holds whenever $f_{1}$ is a linear function and $f_{1},\dots,f_{s-1}$ have a property called totally-homogeneous equidistribution (mod $1$), which we define below.

We say that a function $f:\bN\rightarrow \bR^{s}$ is equidistributed (mod $1$) if for every half-open box $\mathcal{B}=[a_{1},b_{1})\times\cdots [a_{s},b_{s})$, with $0\le a_{i},b_{i}\le 1$,
\[\dens\{n\in \bN: \{f(n)\}\in \mathcal{B}\}=\mu(\mathcal{B}),\]
where $\mu$ denotes the Lebesgue measure and for $\bx=(x_{1},\dots,x_{s})$ we let $\{\bx\}=(\{x_{1}\},\dots,\{x_{s}\})$. We remark that we could equivalently define equidistribution with closed boxes $[a_{1},b_{1}]\times\cdots [a_{s},b_{s}]$ or open boxes $(a_{1},b_{1})\times\cdots (a_{s},b_{s})$ in place of half-open boxes.

We say that $f$ is totally-homogeneously equidistributed (mod $1$) if for every $D\in \bN$ and every $u\in \{0,1,\dots,D-1\}$, we have that $f(Dn+u)/D$ is equidistributed (mod $1$).

If $f_{1},\dots,f_{s-1}$ are totally-homogeneously equidistributed (mod $1$), then in particular, we can show that for any $r\in \bN$,
\[\lim_{N\rightarrow\infty} \dens\{n\in \bN: r\mid \gcd(n,\floor{f_{1}(n)},\dots,\floor{f_{s-1}(n)})\}=\frac{1}{r^{s}}.\]
This follows from the fact that for any $i$, if $n=kr$, then $r\mid \floor{f_{i}(n)}$ if and only if $\{f_{i}(kr)/r\}\in [0,1/r)$.

It follows from the inclusion-exclusion principle that for any $R\ge 1$,
\begin{align*}\dens\{n\in \bN: p\nmid \gcd(n,\floor{f_{1}(n)},\dots,\floor{f_{s-1}(n)})&\text{ for all primes }p\le R\}\\&=\prod_{\substack{p\le R\\ p\text{ prime}}}\left(1-\frac{1}{p^{s}}\right).\end{align*}
For any $x\in \bN$, $x=1$ if and only if $p\nmid x$ for all primes $p$, so if
\[\lim_{R\rightarrow\infty} \overline{\dens}\{n\in \bN:\exists p>R \text{ prime, }p\mid \gcd(n,\floor{f_{1}(n)},\dots,\floor{f_{s-1}(n)})\}=0,\]
then we have
\begin{align*}&\dens \{n\in \bN: \gcd(n,\floor{f_{1}(n)},\dots,\floor{f_{s-1}(n)})=1\}=\\ \lim_{R\rightarrow\infty}&\dens\{n\in \bN: p\nmid \gcd(n,\floor{f_{1}(n)},\dots,\floor{f_{s-1}(n)})\text{ for all primes }p\le R\}\\ &=\prod_{p\text{ prime}}\left(1-\frac{1}{p^{s}}\right)=\zeta(s)^{-1}.\end{align*}
We therefore have the following sufficient conditions for \eqref{main_eq} to hold.

\begin{prop}\label{sufficient_conditions}
Let $f_{1},\dots,f_{s-1}:\bN\rightarrow \bR$ and suppose that:
\begin{enumerate}
\item $(f_{1},\dots,f_{s-1})$ is totally-homogeneously equidistributed (mod $1$) and
\item \[\lim_{R\rightarrow\infty} \overline{\dens}\{n\in \bN:\exists p>R \text{ prime, }p\mid \gcd(n,\floor{f_{1}(n)},\dots,\floor{f_{s-1}(n)})\}=0.\]
\end{enumerate}
Then \eqref{main_eq} holds.
\end{prop}

As shown by Champagne in \cite{c2026}, when $f_{1},\dots,f_{s-1}$ are all polynomials, then we have totally-homogeneous equidistribution (mod $1$) if and only if for every $\bh=(h_{1},\dots,h_{s-1})\in \bZ^{s-1}\setminus\{\mathbf{0}\}$, the polynomial $h_{1}f_{1}(x)+\cdots+h_{s-1}f_{s-1}(x)$ has an irrational coefficient other than the constant coefficient. This in particular holds if the $f_{i}$ are of differing degrees and have irrational leading coefficients, or the same degree and with leading coefficients irrational and linearly independent over $\bQ$.

Champagne then showed that for any $\alpha,\beta\in \bR$ with $\alpha$ irrational,
\[\lim_{R\rightarrow\infty}\overline{\dens}\{n\in \bN: \exists p>R \text{ prime, }p\mid \gcd(n,\floor{\alpha n+\beta})\}=0\]
which gives us \eqref{main_eq} whenever $f_{1}$ is linear and totally-homogeneous equidistribution (mod $1$) holds. Champagne asked whether \eqref{main_eq} holds when $s=2$ and $f_{1}(x)=\alpha x^{2}$ with $\alpha$ irrational and asked more generally whether we have \eqref{main_eq} for any polynomials $f_{1},\dots,f_{s-1}$ which are totally-homogeneously equidistributed (mod $1$). In this paper, we will prove that this is indeed the case.

\begin{thm}\label{main_thm}
Let $s\ge 2$ and $f_{1},\dots,f_{s-1}:\bN\rightarrow \bR$. Suppose that $(f_{1},\dots,f_{s-1})$ is totally-homogeneously equidistributed (mod $1$) and suppose that $f_{1}=\alpha_{d}x^{d}+\cdots+\alpha_{1}x+\alpha_{0}$ for some $d\ge 1$, $\alpha_{0},\dots,\alpha_{d}\in \bR$ with $\alpha_{i}$ irrational for some $1\le i\le d$.

Then
\begin{equation}\label{main_thm_eq}\dens\{n\in \bN:(n,\floor{f_{1}(n)},\dots,\floor{f_{s-1}(n)})=1\}=\frac{1}{\zeta(s)}.\end{equation}

In particular, if $f_{1},\dots,f_{s-1}$ are all polynomials and for every $(h_{1},\dots,h_{s-1})\in \bZ^{s-1}\setminus\{\mathbf{0}\}$, the polynomial $h_{1}f_{1}+\cdots+h_{s-1}f_{s-1}$ has an irrational coefficient other than the constant coefficient, then we have \eqref{main_thm_eq}.
\end{thm}

Note that the $d=1$ case has already been proven in \cite{c2026}, so we need only consider $d\ge 2$. Theorem \ref{main_thm} will follow immediately from Proposition \ref{sufficient_conditions} and the following theorem.

\begin{thm}\label{density_zero}
Let $d\in \bZ$ with $d\ge2$ and $f(x)=\alpha_{d} x^{d}+\alpha_{d-1} x^{d-1}+\cdots+\alpha_{0}$, with $\alpha_{i}\in \bR$ for $0\le i\le d$ and suppose that $\alpha_{i}$ is irrational for some $i\ge 1$. Then
\[\lim_{R\rightarrow\infty}\overline{\dens}\{n\in \bN: \exists p\ge R \: \text{prime},\: p\mid n\:\text{and}\: p\mid \lfloor f(n)\rfloor\}=0.\]
\end{thm}

The subsequent sections will therefore be dedicated to proving Theorem \ref{density_zero}, but first we will say something about a more general question that is still open. Given $f_{1},\dots,f_{s}$ polynomials and $s\ge 2$, we can ask for conditions under which
\begin{equation}\label{main_eq_generalised}\dens\{n\in \bN: \gcd(\floor{f_{1}(n)},\dots,\floor{f_{s}(n)})=1\}=\frac{1}{\zeta(s)}.\end{equation}
In particular, we ask whether totally-homogeneous equidistribution (mod $1$) is still a sufficient condition in this case. Theorem \ref{main_thm} is the special case of the above when $f_{1}(n)=n$. The case of linear $f_{1}(x)=\alpha x+\beta$ also follows from Theorem \ref{density_zero}. This is because we can make a change of variables $n\mapsto \alpha n+\beta$ at the cost of a multiplicative factor of $\max\{|\alpha|,|\alpha|^{-1}\}$ and then apply Theorem \ref{density_zero} with $f=f_{2}(\alpha^{-1}(n-\beta))$.

The fact that, for a prime $p$, $p\mid n$ implies that $n=kp$ for some $k\in \bZ$ turns out to be very helpful. For example, since $\sum_{p\le N} 1/p =\log\log N+O(1)$, which grows very slowly. This means that we only need a small saving on the number of $k$ such that $p\mid \floor{f(kp)}$. If $f_{1}$ has rational coefficients, then for sufficiently large $p$, $p\mid \floor{f_{1}(n)}$ implies that $n=pk+u$ for some $k$, where $u$ is one of at most $\deg f_{1}$ roots of $f_{1}$ (mod $p$). It may therefore be possible to adapt the method of this paper to handle this case as well, although we now have to consider primes up to size $O(N^{\deg f_{1}})$ which causes some difficulties.

In full generality, \eqref{main_eq_generalised} seems to be a difficult problem, and it may require completely different methods. Even showing that $p\mid \floor{f_{1}(n)}$ for at most $O(N/p)$ many $n\in [N]$ is not obvious in this case.

\section*{Notation}
For $x\in \bR$, we write $\floor{x}$ for the largest integer $n\le x$ and write $\{x\}=x-\floor{x}$. We also write $\|x\|$ for the distance from $x$ to the nearest integer. We will use the convention that $\bN=\{1,2,3,\dots\}$.

For $N\in \bN$, we write $[N]=\{n\in \bN: n\le N\}$.

We write $\mu$ for the Lebesgue measure.

Given $c\in \bR$, we write $\bR_{\ge c}=\{x\in \bR: x\ge c\}$ and similarly define $\bR_{>c}$.

Given two functions $f,g$ on variables $x_{1},\dots,x_{m}$, we write $f\ll g$, or equivalently $f=O(g)$, to mean that there is some constant $C$ such that $f\le Cg$, for all relevant values of $x_{1},\dots,x_{m}$. Adding a subscript $f\ll_{y_{1},\dots,y_{n}} g$ or $f=O_{y_{1},\dots,y_{n}}(g)$ indicates that the constant $C$ is allowed to depend on the parameters $y_{1},\dots,y_{n}$.

The variable $p$ will always be taken to be a prime unless specified otherwise.

Whenever we take a sum over the variable $X$ between some parameters $A\le B$, we mean that the sum runs over a dyadic sequence starting at $A$. That is, a sequence of the form $X_{0},X_{1},\dots,X_{m}$, where $X_{0}=A$, $X_{i+1}=2X_{i}$ for all $i$, and $X_{m}\le B$, $2X_{m}$. For example, when we write
\[\sum_{1\le X\le N} X,\]
we mean
\[\sum_{i=0}^{m} X_{i},\]
with $1=X_{0}$, $X_{i+1}=2X_{i}$ for $0\le i\le m-1$ and $X_{m}\le N,2X_{m}$.

\section*{Outline of the proof}

We now outline the idea of the proof of Theorem \ref{density_zero}. We know that at least one of the nonconstant coefficients of $f$ is irrational, and for simplicity, we suppose here that it is the leading coefficient $\alpha_{d}$.

Let $N,R\in \bN$. For a uniform random variable $n$ on $[N]$, we consider the probability that there is some prime $p\ge R$ dividing both $n$ and $\floor{f(n)}$. We wish to prove that $\limsup_{N\rightarrow\infty}$ of this probability tends to $0$ as $R\rightarrow\infty$.

One might hope that the events $p\mid n$ and $p\mid \floor{f(n)}$ are approximately independent, so that the probability that both $p\mid n$ and $p\mid \floor{f(n)}$ occur is $O(p^{-2})$. Observe that if $p\mid n$ and $n=kp$ then $p\mid \floor{f(n)}$ if and only if $mp\le f(kp)< mp+1$, for some integer $m$. This holds if and only if $\{f(kp)/p\}\in [0,1/p)$.

We write $f_{p}(x)=f(px)/p$. The $i$th coefficient of $f_{p}$ is $p^{i-1}\alpha_{i}$, and $p^{d-1}\alpha_{d}$ is irrational, so $f_{p}(k)$ is equidistributed (mod $1$). Therefore, the number of integers $1\le k\le N/p$ such that $\{f_{p}(k)\}\in [0,1/p)$ should be approximately $(N/p)\mu([0,1/p))=Np^{-2}$, which would imply that the probability that $p$ divides both $n$ and $\floor{f(n)}$ is $O(p^{-2})$ as desired. Unfortunately this argument does not quite work, because at the equidistribution step, the rate of convergence could be made arbitrarily slow by making certain choices of $\alpha_{d}$.

To illustrate this, suppose that $\alpha_{d}$ is irrational, but very close to some rational with small denominator and $\alpha_{i}=0$ for $0\le i\le d-1$. Say $|\alpha_{d}-a/q|\le N^{-d}q^{-1}$ with $(a,q)=1$ and $q$ small compared to $N$. Then if $p>2q$ and $1\le k\le N/p$, we have that $|f_{p}(k)-ap^{d-1}k^{d}/q|\le 1/(qp)$. So $\{f_{p}(k)\}\in [0,1/p)$ if and only if $ap^{d-1}k^{d}\equiv 0$ (mod $q$), if and only if $q\mid k^{d}$. There is some $q'\ge q^{1/d}$ such that $q\mid k^{d}$ if and only if $q'\mid k$. So the probability that $p\mid n$ and $p\mid \floor{f(n)}$ will be approximately $1/(q'p)\ge 1/(qp)$, which is considerably larger than the desired $p^{-2}$.

When $\alpha_{d}$ is not too close to a rational with small denomiator, the situation is much nicer. The Erd\H os-Tur\'an inequality tells us that given a sequence $a_{1},\dots,a_{m}$ and an interval $I\subseteq \bR/\bZ$, the number of $1\le i\le m$ such that $\{a_{i}\}\in I$ is approximately $m\mu(I)$ as long as certain exponential sums over $a_{i}$ are small.

Then if $g(x)$ is a polynomial whose leading coefficient is close to a rational $a/q$ with $q$ not too large, then by a certain Weyl sum bound (Lemma \ref{weyl_sum_estimate}), the exponential sum
\[\sum_{1\le k\le M} e(g(k))\]
has a saving of a power of $q$ over the trivial bound $M$.

It can then be shown that if $p^{d-1}\alpha_{d}$ is not too close to any rational with small denominator, then we get a saving of some small power of $N/p$. For these $p$, we find that the probability of $p\mid n$ and $p\mid \floor{f(n)}$ is $O(p^{-2})$ (we prove this in Section \ref{minor_arc_section}). By analogy with the Hardy-Littlewood circle method, we call this the ``minor arc'' contribution.

We refer to the contribution from the remaining $p$ as the ``major arc'' contribution. If $\alpha_{d}$ is not close to a rational with small denominator then the number of such $p$ turns out to be relatively small, and we can get enough saving to prove that the total contribution is $o(1)$. We prove this in Section \ref{major_arc_section}.

In Section \ref{hardest_section}, we complete the proof of Theorem \ref{density_zero} by bounding the major arc contribution in the case where $\alpha_{d}$ is close to a rational with small denominator. Recall that for each $p$, the probability that $p\mid n$ and $p\mid \floor{f(n)}$ is approximately $1/(q'p)$ for some $q^{1/d}\le q'\le q$. Suppose that there is an increasing sequence of integers $N_{1}<N_{2}<\dots$ and approximations $a_{i}/q_{i}$, $(a_{i},q_{i})=1$ with $|\alpha_{d}-a_{i}/q_{i}|\le N_{i}^{-d}q_{i}^{-1}$. We can choose $\alpha_{d}$ in such a way that $q_{i}$ grows arbitrarily slowly compared to $N_{i}$. In particular, we can have $q_{i}\ll \log \log \log N_{i}$. Then
\[\frac{1}{N_{i}}\sum_{2q_{i}\le p\le N_{i}}\frac{1}{q_{i}p}\gg \frac{\log\log N_{i}}{\log\log\log N_{i}},\]
which does not tend to $0$. This may seem hopeless, but in the above bound we haven't taken into consideration the possible overlap between the sets $\{n\in [N]: p\mid n\text{ and } p\mid \floor{f(n)}\}$ for different $p$. We must account for this overlap somehow if we want to prove Theorem \ref{density_zero}.

The idea, and the main innovation of this paper, is to remove a set $E$ from $[N]$ of size $|E|\ll Nq^{-\gamma}$ (for some fixed $\gamma>0$) such that the major-minor arc type argument works when restricted to $[N]\setminus E$. Although the denominator $q$ can grow arbitrarily slowly as $N\rightarrow\infty$, the fact that $\alpha_{d}$ is irrational tells us that we must at least have $q\rightarrow \infty$, so $|E|=o(N)$.

Our construction of such a set $E$ exploits the fact that for a rational polynomial
\[g(x)=(a_{d}x^{d}+\cdots+a_{1}x+a_{0})/q,\]
we have that if $x\equiv t$ (mod $q$) then $\floor{g(x)}=(a_{d}x^{d}+\cdots+a_{1}x+a_{0}-u_{t})/q$, where $u_{t}\in \{0,1,\dots q-1\}$ is such that
\[a_{d}t^{d}+\cdots+a_{0}\equiv u_{t}\: \text{(mod $q$)}.\]
So if $n\equiv t$ (mod $q$) and $p\mid n$, $p\mid \floor{g(n)}$, then $p\mid (a_{0}-u_{t})$. If $p>q$ then there can be at most one possible $0\le u\le q-1$ such that $p\mid (a_{0}-u)$, and therefore $u_{t}=u$. The problem is then reduced to bounding the number of solutions to a polynomial equation (mod $q$).

This is the main idea of the proof, although more work is required to account for the fact that $f$ is not exactly a rational polynomial. We can also apply a major-minor arc analysis to the lower order coefficients of $f$, so that we are left only with the case where all of the coefficients of $f$ except the constant coefficient are close to a rational with small denominator. We then have $f(x)=g(x)+h(x)$ for some rational polynomials $g$ and $h$, where the coefficients of $g$ are rational with small denominator and the coefficients of $h$ are small. Due to the small coefficients, we can break up $[N]$ into suitably nice pieces on which $h$ is almost constant, and then we can treat $f$ like a rational polynomial.
\section{A major-minor arc decomposition}
Given a polynomial $f(x)=\alpha_{d}x^{d}+\cdots+\alpha_{1}x+\alpha_{0}\in \bR[x]$ and $R,N\in \bN$, let
\[\mathcal{S}_{f,R}(N)=\{n\in [N]:\exists p\ge R,\: p\mid \gcd(n,\lfloor f(n)\rfloor)\}.\]
Our aim is to bound the size of $\mathcal{S}_{f,R}(N)$. We will split it into $3$ pieces and bound each piece separately.

Suppose that $X,Y\in \bR$ with $X,Y\ge 1$ and $\eta\in (0,1/2)$. For $q\ge 1$, $a\in \bZ$, and $j\in \bN$, we define
\[\fM_{j}(q,a)=\fM_{j}(X,Y,\eta;q,a)=\{x\in\bR/\bZ:\|x-a/q\|\le Y^{\eta-j}q^{-1}\},\]
\[\fM_{j}=\fM_{j}(X,Y,\eta)=\bigcup_{q\le X^{\eta}}\bigcup_{\substack{0\le a\le q\\ (a,q)=1}} \fM(q,a)\]
and
\[\fm_{j}=\fm_{j}(X,Y,\eta)=(\bR/\bZ)\setminus \fM.\]
Due to analogy with the Hardy-Littlewood circle method, we will call $\fM_{j}$ the major arcs and $\fm_{j}$ the minor arcs.

For $\boldsymbol{\eta}=(\eta_{1},\dots,\eta_{d})\in (0,1/2)^{d}$, let
\[\mathcal{P}_{1}(X,Y,\boldsymbol{\eta})=\{p\in [X,2X]: p\text{ prime, } p^{j-1}\alpha_{j}\in \fM_{j}(X,Y,\eta_{j})\text{ for all }1\le j\le d\}\]
and
\[\mathcal{P}_{2}(X,Y,\boldsymbol{\eta})=\{p\in [X,2X]: p\text{ prime, } p^{j-1}\alpha_{j}\in \fm_{j}(X,Y,\eta_{j})\text{ for some }1\le j\le d\}.\]
We also let, for $i\in \{1,2\}$,
\[\mathcal{D}_{i}(X,Y,\boldsymbol{\eta})=\{n\in [XY]: p\mid n\text{ and } p\mid \floor{f(n)}\text{ for some }p\in \mathcal{P}_{i}(X,Y,\boldsymbol{\eta})\}\]
and
\[D_{i}(X,Y,\boldsymbol{\eta})=\#\mathcal{D}_{i}(X,Y,\boldsymbol{\eta}).\]

Let $\theta\in (0,1/2)$. We define
\[\mathcal{S}_{1}(N,\theta,\boldsymbol{\eta})=\bigcup_{R\le X\le N^{1-\theta}} \mathcal{D}_{1}(X,N/X,\boldsymbol{\eta}),\]
\[\mathcal{S}_{2}(N,\theta,\boldsymbol{\eta})=\bigcup_{R\le X\le N^{1-\theta}} \mathcal{D}_{2}(X,N/X,\boldsymbol{\eta})\]
and
\[\mathcal{S}_{3}(N,\theta)=\{n\in [N]: p\mid n\text{ and } p\mid \floor{f(n)}\text{ for some }p>N^{1-\theta} \}.\]
We find that
\[\mathcal{S}_{f,R}(N)=\mathcal{S}_{1}(N,\theta,\boldsymbol{\eta})\cup \mathcal{S}_{2}(N,\theta,\boldsymbol{\eta})\cup \mathcal{S}_{3}(N,\theta).\]
Note that these three pieces are not necessarily disjoint, as for some $n\in [N]$, there may be multiple primes $p\ge R$ such that $p\mid n$ and $p\mid \floor{f(n)}$.

We will separately bound each of these three pieces to obtain a bound for $\mathcal{S}_{f,R}(N)$. The third piece is the easiest to bound. We have
\begin{align*}\nonumber\sum_{N^{1-\theta}<p\le N} \left(\frac{N}{p}+1\right)= N\log\log N-N\log\log &(N^{1-\theta})+O\left(\frac{N}{\log N}\right)\\ &\ll N\left(\theta+\frac{1}{\log N}\right).\end{align*}
Then it follows that
\[\#\mathcal{S}_{3}(N,\theta)\ll \theta N+\frac{N}{\log N}.\]

It follows that for any $\boldsymbol{\eta}$,
\begin{equation}\label{calS_arc_decomposition}\#\mathcal{S}_{f,R}(N)\ll  \theta N+\frac{N}{\log N}+ \#\mathcal{S}_{1}(N,\theta,\boldsymbol{\eta})+\#\mathcal{S}_{2}(N,\theta,\boldsymbol{\eta}).\end{equation}
By the definition of $D_{1}$ and $D_{2}$, we have, for $i\in \{1,2\}$,
\[\#\mathcal{S}_{i}(N,\theta,\boldsymbol{\eta}) \le \sum_{R\le X\le N^{1-\theta}} D_{i}(X,N/X,\boldsymbol{\eta}),\]
where we recall that the sum over $X$ is over a dyadic sequence. Hence, if we can give a sufficiently good bound on $D_{i}$, then this will give us a bound for $\#\mathcal{S}_{i}(N,\theta,\boldsymbol{\eta})$. We will aim to do this wherever possible, but in Section \ref{hardest_section}, we will not be able to do this.

The main tools we will use will the Erd\H os-Tur\'an inequality and a bound for Weyl sums. We will state these below.

\begin{lemma}\label{weyl_sum_estimate}
Let $N,M,q\in \bZ$ with $M,q\ge 1$. Let $C\ge 1$ and $f(x)=\sum_{k=0}^{d} \alpha_{k}x^{k}$ be a polynomial of degree $d\ge 1$ with $\alpha_{0},\dots,\alpha_{d}\in \bR$ and let $h,a,q\in \bZ$ be such that $q\ge 1$, $(a,q)=1$ and $|q\alpha_{d}-a|\le q^{-1}$. Then for all $\eps>0$,
\begin{equation}\label{weyl}\sum_{x=M}^{M+N} e\left(hf(x)\right)\ll_{\eps,d} N^{1+\eps}\left(\frac{h}{q}+\frac{1}{N}+\frac{(q,h)}{N^{d-1}}+\frac{q}{N^{d}}\right)^{\frac{1}{d(d-1)}}.\end{equation}
The implied constant depends only on $d$ and $\eps$.

Furthermore, if $q\le M$, then under the same assumptions, we have
\begin{equation}\label{weyl_small_q}\sum_{x=M}^{M+N} e\left(hf(x)\right)\ll_{\eps,d} Mh^{\frac{1}{d(d-1)}}q^{\eps-\frac{1}{d(d-1)}}.\end{equation}
\begin{proof}
This is essentially \cite[Lemma 3.2]{st2024}, except that in that paper, the sum goes from $x=1$ to $x=N$. However, since the implied constant depends only on the leading coefficient of $f$, we can shift the sum to the range $M\le x\le M+N$ by replacing $f(x)$ with $f(x-M)$. 

For \eqref{weyl_small_q}, we split the sum into intervals of length between $q/2$ and $q$ and then apply \eqref{weyl} on each interval.
\end{proof}
\end{lemma}

\begin{lemma}[Erd\H os-Tur\'an inequality]\label{erdos_turan}
Let $M,N\in \bN$ and $a_{1},\dots,a_{N}\in \bR$. Then for any interval $I\subseteq \bR/\bZ$, we have
\[\left|\frac{1}{N}\#\{n\in [N]: \{a_{n}\}\in I\}-\mu(I)\right|\ll \frac{1}{M}+\frac{1}{N}\sum_{m=1}^{M} \frac{1}{m}\left|\sum_{n=1}^{N} e(ma_{n})\right|.\]
\begin{proof}
This can be seen for example as a special case of \cite[Theorem 1.21]{dt1997}.
\end{proof}
\end{lemma}

Combining the Weyl sum bound and the Erd\H os-Tur\'an inequality gives us the following corollary.

\begin{cor}\label{weyl_corollary}
Let $N,T\in \bZ$ with $T\ge 1$ and $\eta\in (0,1/2)$. Let $f(x)=\sum_{k=0}^{d} \alpha_{k}x^{k}$ be a polynomial of degree $d\ge 1$ with $\alpha_{0},\dots,\alpha_{d}\in \bR$.

Suppose that there are some $q,a\in \bZ$, $(a,q)=1$, $1\le q\le T^{d-\eta}$ and $|q\alpha_{d}-a|\le q^{-1}$. Let
\[Z=\min\{q,T, T^{d}/q\}.\]
If $d\ge 2$ then for all intervals $I\subseteq \bR/\bZ$ and all $\eps>0$, we have
\[\left|\frac{1}{T}\#\{n\in [N,N+T]: \{f(n)\}\in I\}-\mu(I)\right|\ll_{d,\eta,\eps} Z^{\eps-\frac{1}{2d(d-1)}}.\]
If $d=1$ then we instead have, for all intervals $I\subseteq \bR/\bZ$ and all $\eps>0$,
\[\left|\frac{1}{T}\#\{n\in [N,N+T]: \{f(n)\}\in I\}-\mu(I)\right|\ll Z^{-1/2}.\]
\begin{proof}
First, note that the assumption that $q\le T^{d-\eta}$ means that any bound of $O_{\eps}(T^{\eps})$ or $O_{\eps}(q^{\eps})$ can be replaced with $O_{\eps}(Z^{\eps})$ after changing the value of $\eps>0$.

By the Erd\H os-Tur\'an inequality, for any $M\ge 1$,
\begin{equation}\label{weyl_corollary_erdos_turan}\left|\frac{1}{T}\#\{n\in [N,N+T]: \{f(n)\}\in I\}-\mu(I)\right|\ll \frac{1}{M}+\frac{1}{T}\sum_{m=1}^{M}\frac{1}{m} \left|\sum_{n=N}^{N+T} e(mf(n))\right|.\end{equation}
Suppose first that $d\ge 2$. We let $M=\floor{Z^{1/2}}$. If $q\le T$, then by \eqref{q_saving_bound}, for $1\le m\le M\le q^{1/2}$, we have
\[\sum_{n=N}^{N+T} e(mf(n))\ll_{d,\eps} Tq^{\eps-\frac{1}{2d(d-1)}},\]
which suffices.

Suppose now that $q\ge T$. By Lemma \ref{weyl_sum_estimate}, if $1\le m\le M$ then we have
\[\sum_{n=N}^{N+T} e(mf(n))\ll_{d,\eps} T^{1+\eps} Z^{-\frac{1}{2d(d-1)}}+T^{1+\eps}\left(\frac{(m,q)}{T^{d-1}}\right)^{\frac{1}{d(d-1)}}.\]
The result will follow if we can verify that
\[\sum_{m=1}^{M} \frac{1}{m}T^{1+\eps}\left(\frac{(m,q)}{T^{d-1}}\right)^{\frac{1}{d(d-1)}}\ll_{\eps} q^{\eps}T^{1+\eps-\frac{1}{d}}.\]
To see this, we split the sum according to the value of $(m,q)$. Let $r=(m,q)$ and $m'=m/r$. Then
\[\sum_{m=1}^{M} \frac{1}{m}T^{1+\eps}\left(\frac{(m,q)}{T^{d-1}}\right)^{\frac{1}{d(d-1)}}\le T^{1+\eps}\sum_{r\mid q}\frac{1}{r}\sum_{1\le m'\le q/r}\frac{1}{m'}\left(\frac{r}{T^{d-1}}\right)^{\frac{1}{d(d-1)}}.\]
Now the claimed bound follows upon noting that $\sum_{r\mid q} 1\ll_{\eps} q^{\eps}$.

Now for $d=1$, this is a straightforward application of the Erd\H os-Tur\'an inequality. For $1\le m\le M$, we have
\[\sum_{n=N}^{N+T}e(mf(n))\ll \max\{T,\|m\alpha_{1}\|^{-1}\}.\]
Now, since $|\alpha_{1}-a/q|\le q^{-2}$, if $1\le m \le q/2$ then $\|m\alpha_{1}\|\ge \frac{1}{2}\|am/q\|\ge 1/(2q)$, where for the second inequality we have used the fact that $q\nmid am$. Since $M\le q^{1/2}\le q/2$ unless $q<4$, this tells us that $\|m\alpha_{1}\|^{-1}\ll q$ for all $1\le m\le M$. Then recalling \eqref{weyl_corollary_erdos_turan}, we have
\[\left|\frac{1}{T}\#\{n\in [N,N+T]: \{f(n)\}\in I\}-\mu(I)\right|\ll_{\eps} Z^{-1/2}+(\log Z) q/T\le Z^{-1/2}+Z^{\eps-1}.\]
\end{proof}
\end{cor}

\section{The minor arc contribution}\label{minor_arc_section}
In this section we will bound the size of the minor arc contribution $\mathcal{S}_{2}(N,\theta,\boldsymbol{\eta})$.
\begin{lemma}\label{minor_arc_bound}
Let $d\in \bN$ with $d\ge 2$ and $f(x)=\alpha_{d}x^{d}+\cdots+\alpha_{0}$ with $\alpha_{j}$ irrational for some $1\le j\le d$. Let $X,Y\in \bR$, $N=XY$ and $\theta\in (0,1/2)$ such that $Y\ge N^{\theta}\ge 3$.

Then there exists $\bc=(c_{1},\dots,c_{d})\in \bR_{\ge 1}^{d}$ depending only on $d$ and $\theta$, such that for any $\eta>0$ sufficiently small in terms of $d,\theta$,
\begin{equation}\label{minor_arc_bound_D2}D_{2}(X,Y,\bc \eta)\ll_{d,\eta,\theta} NX^{-\eta/(2d^{2})}.\end{equation}
Consequently, if $R\ge 3$, then
\[\#\mathcal{S}_{2}(N,\theta,\bc \eta)\ll_{d,\eta,\theta} NR^{-\eta/(2d^{2})}.\]
\begin{proof}
Suppose that $p\in [X,2X]$ and let $\mathcal{N}_{p}$ be the number of $k\in [Y]$ such that $p\mid \floor{f(kp)}$. Recall that $p\mid \lfloor f(kp)\rfloor$ if and only if $\{f_{p}(k)\}\in [0,1/p)$, where $f_{p}(k)=f(kp)/p$. The bound \eqref{minor_arc_bound_D2} will follow if we can find some $\bc$ such that for any $1\le j\le d$, if $p^{j-1}\alpha_{j}\in \fm_{j}(X,Y,c_{j}\eta)$ then
\begin{equation}\label{minor_arc_bound_Np}\mathcal{N}_{p}\ll_{d,\eta,\theta}YX^{-\eta/(2d^{2})}.\end{equation}

Suppose first that $p^{d-1}\alpha_{d}\in \fm_{d}(X,Y,\eta)$. The leading coefficient of $f_{p}$ is $p^{d-1}\alpha_{d}$. We will suppose that $\eta<\theta$, so that $X^{\eta}<Y$.

By Dirichlet's approximation theorem, there is some $q\le Y^{d}N^{-\eta}$ and $a\in \bZ$ such that $(a,q)=1$ and $|qp^{d-1}\alpha_{d}-a|\le N^{\eta}Y^{-d}$. Since $p^{d-1}\alpha_{d}\in \fm_{d}$, we must have $q>X^{\eta}$. Then by Corollary \ref{weyl_corollary}, we find that $\mathcal{N}_{p}\ll_{\eps,d}Yp^{-1}+Y X^{\eps-\frac{\eta}{2d(d-1)}}\le YX^{-\eta/(2d^{2})}$ for any $\eps>0$, which implies \eqref{minor_arc_bound_Np}.

We will now use induction on $j$ to prove the claim in the case when $p^{j-1}\alpha_{j}\in \fm_{j}(X,Y,c_{j}\eta)$. Suppose that $1\le j\le d-1$ and that we have found $c_{d},\dots,c_{j+1}\in \bR$ such that $p^{i-1}\alpha_{i}\in \fM_{i}(X,Y,c_{i}\eta)$ for all $j+1\le i\le d$. Suppose that $p^{j-1}\alpha_{j}\in \fm_{j}(X,Y,c_{j}\eta)$ for some $c_{j}\ge 1$ to be chosen later.

There are some $q_{d},\dots,q_{j+1}\in \bN$ such that $q_{i}\le X^{c_{i}\eta}$ and $\|q_{i}p^{i-1}\alpha_{i}\|\le N^{c_{i}\eta}Y^{-i}$ for each $i$. Let $Q=\lcm(q_{d},\dots,q_{j+1})\le X^{(c_{d}+\cdots+c_{j+1})\eta}$ and note that $\|Qp^{i-1}\alpha_{i}\|\le N^{(c_{d}+\cdots+c_{j+1})\eta}Y^{-i}$ for $d\ge i\ge j+1$. Let $a_{d},\dots,a_{j+1}\in \bZ$ be such that $|p^{i-1}\alpha_{i}-a_{i}/Q|$ is minimal for each $j+1\le i\le d$.

Let $C=c_{d}+\cdots+c_{j+1}+1$ and $Y_{j}=YN^{-C\eta}Q^{-1}$. Let $I\subseteq [Y]$ be a set of the form $I=[k_{0},k_{0}+QY_{j})\cap \bZ$ for some $k_{0}\in \bZ$. For each $k\in I$ and $j+1\le i\le d$, we have
\[p^{i-1}\alpha_{i}k^{i}-\frac{a_{i}}{Q}k^{i}=p^{i-1}\alpha_{i}k_{0}^{i}-\frac{a_{i}}{Q}k_{0}^{i}+O(N^{-\eta}).\]
So there is some $c_{I,p}$ such that for any $k\in I$,
\[f_{p}(k)=\frac{a_{d}}{Q}k^{d}+\cdots+\frac{a_{j+1}}{Q}k^{j+1}+p^{j-1}\alpha_{j}k^{j}+\cdots+\alpha_{1}k+c_{I,p}+O(N^{-\eta}).\]
The value of $a_{i}k^{i}/Q$ (mod $1$) depends only on the residue of $k$ (mod $Q$), so if $0\le u\le Q-1$ and $I_{u}=\{k\in \bZ:Qk+u\in I\}$, then $\#I_{u}\le Y_{j}+1$ and there is some interval $J_{I,u,p}\subseteq \bR/\bZ$ with $\mu(J_{I,u,p})\ll p^{-1}+N^{-\eta}\ll X^{-\eta}$ such that if $k\in I_{u}$ and $\{f_{p}(Qk+u)\}\in [0,1/p)$, then
\[\left\{p^{j-1}\alpha_{j}(Qk+u)^{j}+\cdots+\alpha_{1}(Qk+u)\right\}\in J_{I,u,p}.\]
Let $g_{u}(k)=p^{j-1}\alpha_{j}(Qk+u)^{j}+\cdots+\alpha_{1}(Qk+u)$ and note that the leading coefficient is $Q^{j}p^{j-1}\alpha_{j}$. We wish to show that if $c_{j}$ is chosen to be sufficiently large, then $Q^{j}p^{j-1}\alpha_{j}\in \fm_{j}(X,Y_{j},\eta)$.

Since $p^{j-1}\alpha_{j}\in \fm_{j}(X,Y,c_{j}\eta)$, there is no $1\le q\le X^{c_{j}\eta}$ such that $\|qp^{j-1}\alpha_{j}\|\le N^{c_{j}\eta}Y^{-j}$. Now note that $Q\le X^{(c_{d}+\cdots+c_{j+1})\eta}$ and recalling the definition of $Y_{j}$, $Y_{j}\ge Y^{1-C\eta/\theta-(c_{d}+\cdots+c_{j+1})\eta}$. It follows that if we choose $c_{j}$ to be large enough, then $Q^{j}\le X^{(c_{j}-1)\eta}$ and $X^{\eta}Y_{j}^{\eta-j}\le N^{c_{j}\eta}Y^{-j}$, so there is no $1\le q\le X^{\eta}$ such that $\|qQ^{j}p^{j-1}\alpha_{j}\|\le X^{\eta}Y_{j}^{\eta-j}$, and so $Q^{j}p^{j-1}\alpha_{j}\in \fm_{j}(X,Y_{j},\eta)$ as required.

Then by Dirichlet's approximation theorem, we have some $1\le q\le Y_{j}^{j-\eta}X^{-\eta}$, $a\in \bZ$ with $(a,q)=1$ such that $|qQ^{j}p^{j-1}\alpha_{j}-a|\le X^{\eta}Y_{j}^{\eta-j}$. Then $q>X^{\eta}$ and by applying Corollary \ref{weyl_corollary} in the same way as the $j=d$ case, we have
\[\#\{k\in I_{u}: \{g_{u}(k)\}\in J_{I,u,p}\}\ll_{\eps} Y_{j}\mu(J_{I,u,p})+Y_{j}X^{-\eta/(2j^{2})}.\]
Recall that $\mu(J_{I,u,p})\ll X^{-\eta}$. We can cover $[Y]$ with at most $O(Y/(QY_{j}))$ many $I$, and there are exactly $Q$ possible values of $u$, so $\mathcal{N}_{p}\ll_{d,\theta,\eta,\eps} YX^{-\eta/(2j^{2})}$, which implies \eqref{minor_arc_bound_Np}.

Applying this argument inductively on $j$ (starting at $j=d$ and decreasing at each step), we obtain $c_{1},\dots,c_{d}\in \bR_{\ge 1}$ such that either \eqref{minor_arc_bound_Np} holds or there is no $1\le j\le d$ such that $p^{j-1}\alpha_{j}\in \fm_{j}(X,Y,c_{j}\eta)$.

Now recall that $N=XY$, so
\[D_{2}(X,Y,\bc \eta)=\sum_{p\in \mathcal{P}_{2}(X,Y,\bc \eta)} \mathcal{N}_{p}\ll_{d,\eta,\theta} NX^{-\eta/(2d^{2})},\]
as required.

It follows that
\begin{align*}\#\mathcal{S}_{2}(N,\theta,\bc \eta)\le \sum_{R\le X\le N^{1-\theta}} D_{2}(X,Y,\bc \eta) & \ll_{d,\eta,\theta} NR^{- \eta/(2d^{2})}\sum_{k=0}^{\infty} (2^{k})^{-\eta/(2d^{2})}\\ &\ll_{d,\eta,\theta} NR^{-\eta/(2d^{2})}.\end{align*}
\end{proof}
\end{lemma}

\section{The major arc analysis: the easier cases}\label{major_arc_section}
In this section, we will prove some results that allow us to bound the size of $D_{1}(X,Y,\bc\eta)$ is different situations. These will in particular allow us to obtain a satisfactory bound for $\#\mathcal{S}_{1}$ except when the highest degree irrational coefficient of $f$ is close to a rational with very small denominator.

Throughout this section, we will fix some $d\ge 2$ and $f(x)=\alpha_{d}x^{d}+\cdots+\alpha_{1}x+\alpha_{0}\in \bR[x]$. Our first two results will allow us to bound the number of $p$ such that $p^{j-1}\alpha_{j}\in \fM_{j}(X,Y,\bc\eta)$ when $\alpha_{j}$ is not too close to a rational number with small denominator. The first will be more effective when $X$ is small

\begin{lemma}\label{major_arcs_alpha_badapprox}
Let $X,Y\ge 1$, $N=XY$, $\alpha\in \bR$ and $j\in \bN$ with $j\ge 2$. Suppose that $0<\eta<1/2$ and that there is no $m\in \bZ$ such that $1\le m\le X^{3\eta+j-2}$ and $\|m\alpha\|\le N^{3\eta}Y^{-j}$. Then the number of $p\in [X,2X]$ such that $p^{j-1}\alpha \in \fM_{j}(X,Y,\eta)$ is at most
\[O_{j,\eta}(X^{1-\eta}).\]
In particular, suppose that $\bc=(c_{2},\dots,c_{d})\in \bR_{\ge 1}^{d-1}$. If there exists $2\le j\le d$ such that there is no $m\in \bZ$ such that $1\le m\le X^{3c_{j}\eta+j-2}$ and $\|m\alpha_{j}\|\le N^{3c_{j}\eta}Y^{-j}$, then
\[D_{1}(X,Y,\bc \eta)\ll_{d,\eta,\bc} NX^{-\eta}.\]
\begin{proof}
Suppose that $p_{1},p_{2}\in [X,2X]$ such that $p_{1}^{j-1}\alpha\in \fM_{j}(X,Y,\eta;q_{1},a_{1})$ and $p_{2}^{d-1}\alpha\in \fM_{j}(X,Y,\eta,q_{2},a_{2})$ for some $q_{1},q_{2},a_{1},a_{2}$ with $q_{1},q_{2}\le X^{\eta}$ and $(q_{1},a_{1})=(q_{2},a_{2})=1$. Then assuming $N^{\eta}\ge 2$, we have
\[|q_{1}q_{2}(p_{1}^{j-1}-p_{2}^{j-1})\alpha+q_{1}a_{2}-q_{2}a_{1}|\le (q_{1}+q_{2})N^{\eta}Y^{-j}\le N^{3\eta}Y^{-j}.\]
So if $m=q_{1}q_{2}|p_{1}^{j-1}-p_{2}^{j-1}|$ then $\|m\alpha\|\le N^{3\eta}Y^{-j}$. We must therefore have either $m=0$, which implies $p_{1}=p_{2}$; or $m=q_{1}q_{2}|p_{1}^{j-1}-p_{2}^{j-1}|>X^{3\eta+j-2}$. In the second case, $|p_{1}-p_{2}|(j-1)X^{j-2}X^{2\eta}\ge q_{1}q_{2}|p_{1}^{j-1}-p_{2}^{j-1}|> X^{3\eta+j-2}$, so $|p_{1}-p_{2}|>X^{\eta}$ and there are at most $O(X^{1-\eta})$ values of $p\in [X,2X]$ such that $p^{j-1}\alpha\in \fM_{j}$.

The second part follows immediately from the definition of $D_{1}(X,Y,\bc \eta)$.
\end{proof}
\end{lemma}

Note that when $X^{3\eta+j-2}\ge N^{3\eta}Y^{-j}$, this is trivial by Dirichlet's approximation theorem. We will therefore need a different bound that is applicable when $X$ is larger, which we do in the next lemma. Note also that the assumption that $j\ge 2$ is essential in the above. When $j=1$, we have $p^{j-1}\alpha=\alpha$, which does not depend on $p$. We find that if $\alpha_{1}\in \fm_{1}(X,Y,c_{1}\eta)$ then $D_{1}(X,Y,\bc\eta)=0$, so we can always assume that $\alpha_{1}\in \fM_{1}(X,Y,c_{1}\eta)$.

\begin{lemma}\label{alpha_goodapprox_bound_X_large}
Let $N\in \bN$, $1\le X\le N$, $Y=N/X$, $\theta\in (0,1/2)$ such that $Y\ge N^{\theta}$ and let $\bc=(c_{1},\dots,c_{d})\in \bR_{\ge 1}^{d}$. Let $\bc'=(c'_{2},\dots,c'_{d})$ be defined by $c'_{j}=3(j-1)^{2}c_{j}$ for $j\ge 2$. Let $\eta\in (0,1/2)$ be sufficiently small in terms of $d,\theta$.

Suppose that there exists $2\le j\le d$ such that there is no $m\in \bZ$, $1\le m\le X^{c'_{j}\eta}$ such that $\|m\alpha_{j}\|\le X^{c'_{j}\eta+1-j}$. Then
\[D_{1}(X,Y,\bc \eta)\ll_{d,\eta} NX^{-\eta}.\]
\begin{proof}
Consider some minor arc $\fM_{j}(q,a)$, which is an interval in $\bR/\bZ$ of length $2N^{c_{j}\eta}Y^{-j}q^{-1}$. By Dirichlet's approximation theorem, there is some $1\le m\le X^{j-1-c'_{j}\eta}$ and $a\in \bZ$ with $(a,m)=1$ such that $|m\alpha_{j}-a|\le X^{c'_{j}\eta+1-j}$. By assumption, we must have $m>X^{c'_{j}\eta}$. By Corollary \ref{weyl_corollary} applied to the polynomial $\alpha_{j}x^{j}$, with $N=T=X$,
\[\left|\frac{1}{X}\#\{X\le r\le 2X: r^{j-1}\alpha_{j}\in \fM_{j}(q,a)\}-\frac{2N^{\eta}}{qY^{j}}\right|\ll X^{-c'_{j}\eta/(2(j-1)^{2})}.\]
Noting that $c'_{j}/(2(j-1)^{2})=3c_{j}$, this gives
\[\#\{X\le r\le 2X: r^{d-1}\alpha\in \fM(q,a)\}\ll_{\eps} X^{1-3c_{j}\eta}.\]
Now note that there are at most $X^{2c_{j}\eta}$ choices for $q$ and $a$, and at most $Y$ choices of $n\in [N]$ with $p\mid n$ for each $p$.
\end{proof}
\end{lemma}

Suppose that $\alpha_{j}$ is irrational and $\alpha_{d},\dots,\alpha_{j+1}$ are rational with common denominator $Q$, say. By restricting to a congruence class (mod $Q$), we can reduce to the case of a degree $j$ polynomial, as we did in the proof of Lemma \ref{minor_arc_bound}. Suppose that $\alpha_{j}\approx a/q$ with $(q,a)=1$. We wish to prove a bound on $D_{1}(X,Y,\bc\eta)$ which saves a power of $q$ over the trivial bound $XY/\log X$. We will show this in Lemma \ref{major_arc_bound_combined_large_X}, where we will deduce it from a more general result: Lemma \ref{q_saving_bound}. We will also use this lemma in the next section.

Note that if $g(x)$ is a polynomial of degree $j$ and the leading coefficient lies in some major arc $\fM_{j}(r,b)$, then by Corollary \ref{weyl_corollary}, the number of $k\in [Y]$ such that $\{g(x)\}$ lies in some interval $J$ is at most $O(Y\mu(J)+Yr^{\eps-2^{1-j}})$. In Lemma \ref{q_saving_bound}, for each prime $X\le p\le 2X$, we will have a polynomial $g_{p}$ with leading coefficient $p^{j-1}\alpha$, $\alpha\approx a/q$ and we wish to bound the number of pairs $(p,k)$ with $\{g_{p}(k)\}\in J$ and $p^{j-1}\alpha\in \fM_{j}$. By the above, the major arcs $\fM_{j}(r,b)$ with $r\gg q$ are not a problem, so we just need to bound the number of $p$ lying in the major arcs $\fM_{j}(r,b)$ with $r$ small. For this, we need the following two lemmas.

\begin{lemma}\label{mod_kth_roots}
Let $q,d\in \bN$ and $y\in \bZ/q\bZ$. Then for all $\eps>0$, the number of solutions $x\in (\bZ/q\bZ)^{\times}$ to the equation
\[x^{k}\equiv y\: \text{(mod $q$)}\]
is at most $O_{d,\eps}(q^{1/2+\eps})$.
\begin{proof}
Using the identity
\[\sum_{a=0}^{q-1}e(ax/q)=\begin{cases}q,\quad \text{if }q\mid x\\ 0,\quad \text{otherwise},\end{cases}\]
we can write the number of solutions as
\begin{equation}\label{mod_kth_roots_eq1}\sum_{x\in (\bZ/q\bZ)^{\times}} \frac{1}{q} \sum_{y=0}^{q-1} e(yx^{k}/q)=\frac{1}{q}\sum_{y=0}^{q-1}\left(\sum_{x\in (\bZ/q\bZ)^{\times}}e(yx^{k}/q)\right).\end{equation}
We have the bound,
\begin{equation}\label{mod_kth_roots_eq2}\left|\sum_{x\in (\bZ/q\bZ)^{\times}}e(yx^{k}/q)\right|\ll_{\eps} q^{1/2+\eps}(q,y)^{1/2},\end{equation}
(see, for example, \cite[Lemma 2.3]{h2026}). For $u\mid q$, there are less than $q/u$ values of $0\le y\le q-1$ for which $(q,y)=u$, so for any $\lambda\le 1$,
\begin{equation}\label{gcd_q_sum}\sum_{y=0}^{q-1} (q,y)^{\lambda}\le \sum_{u\mid q} (q/u)u^{\lambda}=q\sum_{u\mid q} u^{\lambda-1}\le \tau(q)q\ll_{\eps} q^{1+\eps},\end{equation}
where $\tau(q)$ is the number of divisors of $q$.

The result now follows from combining \eqref{mod_kth_roots_eq1}, \eqref{mod_kth_roots_eq2} and \eqref{gcd_q_sum}.
\end{proof}
\end{lemma}

\begin{lemma}[Siegel-Walfisz in short intervals]\label{siegel_walfisz}
For $x,q\in \bN$ and $a\in \bZ$ with $(q,a)=1$, let $\pi_{q,a}(x)$ be the number of primes $p\le x$ such that $p\equiv a$ (mod $q$). Let $7/12<\lambda< 1$, $A, B\ge 1$ be fixed constants.

Then if $1\le q\le (\log x)^{B}$ and $y\le x^{\lambda}$, we have
\[\pi_{q,a}(x+y)-\pi_{q,a}(x)=\frac{y}{\varphi(q)\log x}+O_{A,B,\lambda} \left(\frac{x^{\lambda}}{\varphi(q)(\log x)^{A}}\right),\]
where $\varphi$ is Euler's totient function.
\begin{proof}
Let
\[\psi_{q,a}(x)=\sum_{\substack{n\le x\\ n\equiv a\text{ (mod $q$)}}}\Lambda(n),\]
where
\[\Lambda(n)=\begin{cases}\log p\quad \text{if $n=p^{k}$ for some prime $p$ and $k\in \bN$},\\ 0\quad \text{otherwise}\end{cases}\]
is the Von-Mangoldt function.

Then it follows from \cite[Theorem 2]{z1989} that for all $0\le y\le x^{\lambda}$,
\[\psi_{q,a}(x+y)-\psi_{q,a}(x)=\frac{y}{\varphi(q)}+O_{A,B,\lambda} \left(x^{\lambda}(\log x)^{-A-B})\right).\]
By a standard argument, the contribution to $\psi_{q,a}(x+y)-\psi_{q,a}(x)$ coming those $n=p^{k}$ with $p$ a prime and $k\ge 2$ is at most $\ll ((x+y)^{1/2}-x^{1/2}+1)\log x\ll (y x^{-1/2}+1)\log x\ll x^{1/2}\log x$. For each prime $x\le p\le x+y$, we have $\log p=\log x+O(x^{\lambda-1})$ and it can be deduced that
\[\psi_{q,a}(x+y)-\psi_{q,a}(x)=(\log x)(\pi_{q,a}(x+y)-\pi_{q,a}(x))+O\left(yx^{\lambda-1}+x^{1/2}\log x\right).\]
The result follows, noting that $yx^{\lambda-1}\le x^{2\lambda-1}\ll_{\lambda,A} x^{\lambda}(\log x)^{-A}$.
\end{proof}
\end{lemma}

\begin{lemma}\label{q_saving_bound}
Let $N\in \bN$, $1\le X\le N$, $Y=N/X$ and $j\in \bN$, $j\ge 2$, $\theta\in (0,1/2)$ and $0<\kappa\le j-2$ such that $Y>\max\{X^{2\kappa},N^{\theta}\}$. Consider some fixed $\eta>0$ sufficiently small in terms of $d,\theta$.

Suppose that $\alpha\in \bR$, $q\in \bN$ and $a\in \bZ$ such that $(a,q)=1$, $q\le X^{\kappa}$ and $|q\alpha-a|\le X^{-3\eta+4/3-j-\kappa}$. For each $p\in [X,2X]$, suppose that $g_{p,u}(x)\in \bR[x]$ is a polynomial of degree $j$ with leading coefficient $p^{j-1}\alpha$. For $J\subseteq \bR/\bZ$ an interval, let $\mathcal{N}(J)$ be the number of pairs $(p,k)$ with $p\in [X,2X]$ prime and $k\in [Y]$ such that $p^{j-1}\alpha\in \fM_{j}(X,Y,\eta)$ and $\{g_{p,u}(k)\}\in J$. Then for any $J$ and $\eps>0$,
\begin{equation}\label{q_saving_bound_eq1}\mathcal{N}(J)\ll_{j,\theta,\eta} NX^{-1/2}+N q^{\eps-\frac{1}{2j(j-1)}}(\log X)^{-1}+N\mu(J)(\log X)^{-1}.\end{equation}
\begin{proof}
We divide $[X,2X]$ into $O(X^{1/3})$ many intervals $I$ of length at most $X^{2/3}$. Consider some $I$. We can suppose that there is $r\le X^{\eta}$ and $0\le b\le r-1$ with $(b,r)=1$ and some prime $p_{0}\in I$ such that $\{p_{0}^{j-1}\alpha\}\in \fM_{j}(r,b)$, because otherwise the contribution to $\mathcal{N}(J)$ from those $p\in I$ is $0$. If $0\le t\le q-1$ and $p\in I$ are such that $(p^{j-1}-p_{0}^{j-1})a\equiv t$ (mod $q$), then $|p^{j-1}-p_{0}^{j-1}|\le (j-1)X^{j-4/3}$, so
\[\left\|p^{j-1}\alpha-\frac{t}{q}-\frac{b}{r}\right\|\le (j-1)X^{j-4/3}|\alpha-a/q|+r^{-1}N^{\eta}Y^{-j}\ll X^{-\kappa-3\eta},\]
where for the final inequality we have used the assumption that $Y>X^{2\kappa}$.

This is at most $X^{-\kappa-2\eta}/2$ when $X$ is sufficiently large in terms of $\eta$, and when $X\ll_{\eta} 1$, we have $NX^{-1}\gg_{\eta} N$, so \eqref{q_saving_bound_eq1} holds trivially. We claim that if $s\le X^{\eta}$, $0\le c\le s-1$ such that $(c,s)=1$ and $\{p^{j-1}\alpha\}\in \fM_{j}(s,c)$, then $c/s\equiv b/r+t/q$ (mod $1$).

Indeed, we have $\|p^{j-1}\alpha-c/s\|\le N^{\eta}Y^{-j}<X^{-2\eta-\kappa}/2$. So if $c/s\not\equiv t/q+b/r$ (mod $1$), then $X^{-2\eta-\kappa}\le (qrs)^{-1}\le \|c/s-t/q-b/r\|<X^{-2\eta-\kappa}$, which is a contradiction.

Letting $r'=r/(q,r)$, there is some $b'\in \bZ$ such that $t/q+b/r=(tr'+b')/(qr')$. If $u=(tr'+b',q)$, then by Corollary \ref{weyl_corollary}, the number of $k\in [1,Y]$ such that $\{f_{p}(k)\}\in J$ is at most $O_{\eps}(Y\mu(J)+Y(q/u)^{\eps-1/(2j(j-1))})$ for any $\eps>0$. Furthermore, since $(p^{j-1}-p_{0}^{j-1})a\equiv t$ (mod $q$) and $(a,q)=1$, if $(tr'+b',q)=u$, then $p^{j-1}\equiv y_{I,u}$  (mod $u$) for some $y_{I,u}$ depending only on $I$ and $u$. If $p\nmid u$ then $p$ is invertible (mod $u$) and by Lemma \ref{mod_kth_roots}, $p$ must lie in one of at most $\ll_{j,\eps} u^{1/2+\eps}$ congruence classes (mod $u$). There are also at most $\ll_{\eps} u^{\eps}$ many $p$ such that $p\mid u$.

Fix some $\eps>0$. If $u>(\log X)^{1/\eps}$, then the number of $p\in I$ such that $(tr'+b',q)=u$ is therefore at most $\ll_{j,\eps}u^{\eps}+X^{2/3}u^{\eps-1/2}\le u^{\eps}+X^{2/3}u^{2\eps-1/2}(\log X)^{-1}$. If $u\le (\log X)^{1/\eps}$, then by Lemma \ref{siegel_walfisz} with $\lambda=2/3$, $A=1$ and $B=1/\eps$, the number of $p\in I$ such that $(tr'+b',q)=u$ is again at most $\ll_{j,\eps}u^{\eps}+X^{2/3}u^{2\eps-1/2}(\log X)^{-1}$. In either case, recalling that there are $O(X^{1/3})$ intervals $I$, we find that
\begin{align*}\mathcal{N}(J)&\ll_{j,\eps,\theta,\eta} \sum_{X\le p\le 2X} Y\mu(J)\\+&\sum_{u\mid q}\left(X^{1/3}Yu^{\eps}(q/u)^{\eps-\frac{1}{2j(j-1)}}+XYu^{\frac{1}{2j(j-1)}+2\eps-1/2}q^{\eps-\frac{1}{2j(j-1)}}(\log X)^{-1}\right).\end{align*}
Now, $\sum_{u\mid q}u^{\lambda}\ll_{\eps} q^{\eps}$ for any $\lambda\le 0$ and $\eps>0$. Since $q\le X^{\kappa}$, if $\eps$ is sufficiently small in terms of $\kappa$ then $X^{1/3}Y\sum_{u\mid q}u^{\eps}(q/u)^{\eps-1/(2j(j-1))}\ll_{\eps} X^{1/2}Y=NX^{-1/2}$. Finally, $\sum_{X\le p\le 2X} Y\mu(J)\ll XY(\log X)^{-1}\mu(J)=N\mu(J)(\log X)^{-1}$. The result follows.
\end{proof}
\end{lemma}

We now deduce a bound for $D_{1}(X,Y,\bc\eta)$.

\begin{lemma}\label{major_arc_bound_combined_large_X}
Let $\theta\in (0,1/d)$. Suppose that there exist $1\le j\le d$ such that $\alpha_{j}$ is irrational and suppose that $j$ is maximal with this property. Let $q\ge 1$ be the smallest positive integer such that $\|q\alpha_{j}\|\le N^{-\theta/2}$.

Let $N\in \bN$, $1\le X\le N$, $Y=N/X$ with $Y\ge N^{\theta}$, $\bc=(c_{1},\dots,c_{d})\in \bR_{\ge 1}^{d}$ and let $\eta>0$ be sufficiently small in terms of $d,\theta$. Then for $j\ge 2$,
\begin{equation}\label{major_arc_bound_combined_large_X_eq1}D_{1}(X,Y,\bc \eta)\ll_{f,\eta,\bc} NX^{-\eta}+\frac{Nq^{-1/(2j^{2})}}{\log X},\end{equation}
and when $j=1$,
\begin{equation}\label{major_arc_bound_combined_large_X_eq2}D_{1}(X,Y,\bc \eta)\ll_{f,\eta,\bc} NX^{-\eta}+\frac{Nq^{-1/2}(\log q)}{\log X}.\end{equation}
\begin{proof}
Let $Q,a_{d},\dots,a_{j+1}\in \bZ$ be such that $\alpha_{i}=a_{i}/Q$ for $j+1\le i\le d$. For $0\le u\le Q-1$ and $p\in [X,2X]$ let $c_{u}=(a_{d}p^{d-1}u^{d}+\cdots+a_{j+1}p^{j}u^{j+1})/Q$ and $g_{p,u}(x)=\alpha_{j}p^{j-1}(Qx+u)^{j}+\cdots+\alpha_{1}(Qx+u)+\alpha_{0}p^{-1}+c_{u}$. Then for any $p\in [X,2X]$ and $k\in \bZ$, we have
\[\{f_{p}(Qk+u)\}=\{g_{p,u}(k)\}.\]
So $p\mid \floor{f(p(Qk+u))}$ if and only if $\{g_{p,u}(k)\}\in [0,1/p)$. Let $Y'=\ceil{Y/Q}$, $I_{u}=\{k\in \bZ: Qk+u\in [Y]\}$ and note that $\#I_{u}\le Y'$. Since $Q$ depends only on $f$, we can suppose that it is sufficiently small in terms of $N$ and so $Y'\ge N^{\theta-\eps}$. Suppose first that $j\ge 2$.

By the definition of $Y'$, for $p\in [X,2X]$, if $p^{j-1}\alpha_{j}\in \fM_{j}(X,Y,c_{j}\eta)$ then $p^{j-1}Q^{j}\alpha_{j}\in \fM_{j}(X,Y',c_{j}\eta+\eps)$, assuming that $Y'$ is large enough in terms of $\eps$. The $+\eps$ is needed here because $Y'$ may be slighlty larger than $Y/Q$.

We may suppose that there is $1\le m\le X^{3c_{j}\eta+j-2}$ such that $\|m\alpha_{j}\|\le N^{3c_{j}\eta}Y^{-j}$, because otherwise by Lemma \ref{major_arcs_alpha_badapprox} we have $D_{1}(X,Y,\bc\eta)\ll_{d,\eta} NX^{-\eta}$. If $Y'>X^{2(j-2)}$ then by Lemma \ref{q_saving_bound} with $\kappa=j-2$, the number of $k\in I_{u}$ and $p\in [X,2X]$ such that $p^{j-1}Q^{j}\alpha_{j}\in \fM_{j}(X,Y',c_{j}\eta+\eps)$ and $\{g_{p,u}(k)\}\in [0,1/p)$ is at most
\[\ll_{f,\theta,\bc,\eta}NQ^{-1}X^{-1/2}+NQ^{-1}m^{\eps-1/(2j(j-1))}(\log X)^{-1}.\]
We have that $mQ^{j}\ge q$ because $\|mQ^{j}\alpha_{j}\|\le N^{3\eta}(Y')^{-j}\le N^{-\theta/2}$. So $m\gg_{\eps} q^{1-\eps}$ and \eqref{major_arc_bound_combined_large_X_eq1} follows in this case.

If $Y'\le X^{2(j-2)}$ then $X\ge N^{1/(2j-1)}Q^{-1}$ and by Lemma \ref{alpha_goodapprox_bound_X_large}, we may suppose that there is some $1\le m\le X^{c'_{j}\eta}$ such that $\|m\alpha_{j}\|\le X^{c'_{j}\eta+1-j}$, because otherwise $D_{1}(X,Y,\bc\eta)\ll_{d,\eta} NX^{-\eta}$. Then $\|mQ^{j}\alpha_{j}\|\le X^{c'_{j}\eta+\eps+1-j}\le N^{-\theta/2}$, so $q\le mQ^{j}$ so $m\gg_{\eps} q^{1-\eps}$. By supposing that $\eta$ is sufficiently small in terms of $\theta$, we have $Y'\ge X^{2c'_{j}\eta+\eps}$, so by Lemma \ref{q_saving_bound} with $\kappa=c'_{j}\eta+\eps$, the number of $k\in I_{u}$ and $p\in [X,2X]$ such that $p^{j-1}Q^{j}\alpha_{j}\in \fM_{j}(X,Y',c_{j}\eta+\eps)$ and $\{g_{p,u}(k)\}\in [0,1/p)$ is at most
\[\ll_{f,\theta,\eta,\bc,\eps}NQ^{-1}X^{-1/2}+NQ^{-1}m^{\eps-1/(2j(j-1))}(\log X)^{-1}.\]
Recalling that $m\gg_{\eps} q^{1-\eps}$, we again have \eqref{major_arc_bound_combined_large_X_eq1}.

Now suppose that $j=1$. In this case, for $1\le m< q/Q$, we have $\|mQ\alpha_{1}\|>N^{-\theta/2}$, so
\[\left|\sum_{k\in I_{u}} e(mg_{p,u}(k))\right|\le \|mQ\alpha_{1}\|^{-1}<N^{\theta/2}.\]
and by the Erd\H os-Tur\'an inequality with $M=q/Q$, the number of $k\in I_{u}$ such that $\{g_{p,u}(k)\}\in [0,1/p)$ is at most 
\[\ll Y'p^{-1}+Y'Qq^{-1}+(\log q)N^{\theta/2}.\]
By Dirichlet's approximation theorem and the definition of $q$, we find that $q\le N^{\theta/2}$, so $(\log q)N^{\theta/2}\ll (\log N)Y'N^{-\theta/2}$. Summing over the $O(X/(\log X))$ choices of $p$ and $Q$ choices of $u$ gives \eqref{major_arc_bound_combined_large_X_eq2}.
\end{proof}
\end{lemma}

Summing the bound \eqref{major_arc_bound_combined_large_X_eq1} over $X$ gives the following.

 \begin{cor}\label{major_arcs_alpha_goodapprox_large_X}
Suppose that $N\in \bN$ and let $\theta\in (0,1/2)$, $\lambda\in (0,1-\theta)$, $\bc\in \bR_{\ge 1}^{d}$ and $\eta>0$ with $\eta$ sufficiently small in terms of $d,\theta$ and $\bc$.

Suppose that there exists $1\le j\le d$ such that $\alpha_{j}$ is irrational and $\alpha_{i}$ is rational for all $i>j$. Let $q\ge 1$ be the smallest positive integer such that $\|q\alpha_{j}\|\le N^{-\theta/2}$. Then
\begin{equation}\label{large_X_bound_major_arcs}\sum_{N^{\lambda}\le X\le N^{1-\theta}}D_{1}(X,N/X,\bc\eta)\ll_{f,\eta,\bc,\theta,\lambda} Nq^{-1/(2d^{2})}+N^{1-\eta\lambda}\end{equation}
and if $1\le R\le N^{\theta/2}$, then
\begin{equation}\label{large_q_bound}\sum_{R\le X\le N^{1-\theta}}D_{1}(X,N/X,\bc\eta)\ll_{f,\eta,\bc,\theta} Nq^{-1/(2d^{2})}\log\log N+NR^{-\eta},\end{equation}
where both sums are over dyadic sequences.
\begin{proof}
Note that for $1\le A\le 2B$, we have the following bounds for sums over a dyadic sequence:
\[\sum_{A\le X\le B} X^{-1}\ll A^{-1}\]
and
\[\sum_{A\le X\le B} (\log X)^{-1}\ll \sum_{\log_{2}A\le j\le \log_{2} B}\frac{1}{j}\ll \log\log B-\log\log A+1,\]
where we write $\log_{2} x=(\log x)/(\log 2)$. When $A=N^{\lambda}$ and $B=N^{1-\theta}$, we have $\log\log N^{1-\theta}-\log\log N^{\lambda}=\log (1-\theta)-\log \lambda\ll_{\theta,\lambda} 1$ and when $B=N^{1-\theta}$, $A=R$, we have $\log\log N^{1-\theta}-\log\log R\ll \log\log N$. The result then follows from Lemma \ref{major_arc_bound_combined_large_X}.
\end{proof}
\end{cor}

The first bound is always $o(N)$ as $N\rightarrow\infty$ with $\lambda,\eta,f$ fixed, so this completely handles the contribution to $\#\mathcal{S}_{1}$ coming from $X\ge N^{\lambda}$. However, \eqref{large_q_bound} is only $o(N)$ when $q(\log\log N)^{-2d^{2}}\rightarrow\infty$ as $N\rightarrow\infty$. This is already sufficient to prove Theorem \ref{density_zero} in the case where $\alpha_{j}$ is of finite exponential type. In the next section, we will deal with the remaining cases so that we can prove Theorem \ref{density_zero} in full generality.

\section{The major arc analysis: the hardest case}\label{hardest_section}
In this section we will bound $\#\mathcal{S}_{1}(N,\theta,\bc\eta)$ in the general case. The first step is to show that $D_{1}(X,Y,\bc\eta)$ is small unless all of the nonconstant coefficients of $f$ are very close to rationals with small denominators. Our argument only works when $X\le Y^{1/(2d)}$, but this is not a problem because of \eqref{large_X_bound_major_arcs} bounds the contribution from the remaining $X$.

As in the previous section, we will fix $d\ge 2$ and $f(x)=\alpha_{d}x^{d}+\cdots+\alpha_{0}\in \bR[x]$ with $\alpha_{j}$ irrational for some $1\le j\le d$. 
\begin{lemma}\label{lower_degree_coefficients}
Suppose that $X,Y\ge 2$, $XY=N$ with $Y\ge X^{2d}$, and $\bc=(c_{1},\dots,c_{d})\in \bR_{\ge 1}^{d}$. Let $\eta \in (0,1/2)$ be sufficiently small in terms of $d, \bc$. Then there exists $\sigma=\sigma(d,\bc)>0$ and $\gamma=\gamma(d,\bc)>0$, depending only on $d,\bc$ such that either
\begin{equation}\label{lower_degree_coefficients_eq1}D_{1}(X,Y,\bc \eta)\ll_{d,\bc,\eta} NX^{-\eta}\end{equation}
or there exists $Q\in \bN$, $Q\le X^{\sigma\eta}$ such that $\|Q\alpha_{i}\|\le N^{\sigma\eta}Y^{-i}$ for all $1\le i\le d$ and
\begin{equation}\label{lower_degree_coefficients_eq3}D_{1}(X,Y,\bc\eta)\ll_{d,\bc,\eta} NX^{-\eta}+NQ^{-\gamma}(\log X)^{-1}.\end{equation}
\begin{proof}
By Lemma \ref{major_arcs_alpha_badapprox}, either \eqref{lower_degree_coefficients_eq1} holds, or there is some integer $1\le q\le X^{3c_{d}\eta+d-2}$ such that $\|q\alpha_{d}\|\le N^{3c_{d}\eta}Y^{-d}$. Recall the notation $f_{p}(k)=f(kp)/p$. Then $p\mid \lfloor f(kp)\rfloor$ if and only if $\{f_{p}(k)\}\in [0,1/p)$, so by Lemma \ref{q_saving_bound} with $\kappa=d-2+c_{d}\eta$, $J=[0,1/p)$ and $g_{p}=f_{p}$, noting that $Y\ge X^{2d}>X^{2\kappa}$ if $\eta$ is sufficiently small, we have the bound
\[D_{1}(X,Y,\bc\eta)\ll_{d,\eps,\bc,\eta} NX^{-1/2}+N q^{-1/(2d^{2})}(\log X)^{-1}.\]
We then have \eqref{lower_degree_coefficients_eq1} unless $q\le X^{2d^{2}\eta}$.

Suppose that $1\le j\le d-1$ and $Q,a_{d},\dots,a_{j+1}\in \bZ$, $\sigma,\gamma>0$ are such that $|Q\alpha_{i}-a_{i}|\le N^{\sigma \eta}Y^{-i}$ for $j+1\le i\le d$, $(Q,a_{d},\dots,a_{j+1})=1$, $Q\le X^{\sigma \eta}$ and
\[D_{1}(X,Y,\bc\eta)\ll_{d,\eps,\bc,\eta} NX^{-\eta}+NQ^{-\gamma}(\log X)^{-1}.\]
If we can prove that there exist $\sigma',\gamma'>0$ depending only on $d,\bc,\sigma,\gamma$ such that either \eqref{lower_degree_coefficients_eq1} holds or there exists and $Q'\le X^{\sigma'\eta}$ such that $\|Q'\alpha_{i}\|\le N^{\sigma' \eta}Y^{-i}$ for $j\le i\le d$ and
\begin{equation}\label{lower_degree_coefficients_eq5}D_{1}(X,Y,\bc\eta)\ll_{d,\eps,\bc,\eta} NX^{-\eta}+N(Q')^{-\gamma'}(\log X)^{-1},\end{equation}
then the desired result will follow by induction on $j$ (with the induction starting at $j=d$ and then decreasing at each step).

From now on, suppose that $Q$ is as above and $1\le j\le d-1$. Let $Y_{j}=YN^{-(\sigma+1) \eta}Q^{-1}+1$ and note that for any $p$, if $p^{j-1}\alpha_{j}\in \fM_{j}(X,Y,c_{j}\eta)$, then $p^{j-1}Q^{j}\alpha_{j}\in \fM_{j}(X,Y_{j},c_{j}\eta)$. We restrict our attention to an interval $I=[k_{0},k_{0}+Y_{j}Q]\cap\bZ$ for some $k_{0}\in \bZ$ with $0\le k_{0}$, $k_{0}+Y_{j}Q\le NX^{-1}$ and let $\mathcal{N}_{I}$ be the number of pairs $(p,k)$ with $X\le p\le 2X$ prime, $p^{j-1}\alpha_{j}\in \fM_{j}(X,Y,c_{j}\eta)$ and $k\in I$ such that $p\mid \floor{f(kp)}$. Note that we can cover $[Y]$ by at most $O(Y/(QY_{j}))$ many of these $I$.

Consider some $p\in [X,2X]$. For $k\in I$, we have $p\mid \floor{f(kp)}$ if and only if $\{f_{p}(k)\}\in [0,1/p)$. For $k\in I$, and $i\ge j+1$, we have $k^{i}-k_{0}^{i}\ll k^{i-1}|k-k_{0}|\ll Y^{i-1}Y_{j}Q\ll Y^{i}N^{-(\sigma+1)\eta}$. For each $j+1\le i\le d$, recalling that $|Q\alpha_{i}-a_{i}|\le N^{\sigma\eta}Y^{-i}$, we find that
\[\left|\left(\alpha_{i}-\frac{a_{i}p^{i-1}}{Q}\right)k^{i}-\left(\alpha_{i}-\frac{a_{i}p^{i-1}}{Q}\right)k_{0}^{i}\right|\ll  N^{-\eta}\ll p^{-\eta}.\]
It follows that there is some interval $J_{p,I}$ with $\mu(J_{p,I})\ll p^{-\eta}$ such that if $\{f_{p}(k)\}\in [0,1/p)$ then
\[\left\{\frac{a_{d}p^{d-1}}{Q}k^{d}+\cdots+\frac{a_{j+1}p^{j}}{Q}k^{j+1}+\alpha_{j}p^{j-1}k^{j}+\cdots+\alpha_{0}p^{-1}\right\}\in J_{p,I}.\]
For each $d\ge i\ge j+1$, the value of $a_{i}p^{i-1}k^{i}/Q$ (mod $1$) depends only on the residue of $k$ (mod $Q$). For each $0\le u\le Q-1$, let $I_{u}=\{k\in \bZ: Qk+u\in I\}$ and note that $|I_{u}|\le Y_{j}$. Then there exists an interval $J_{I,u,p}\subseteq \bR/\bZ$ which is a translate of $J_{I,p}$, such that if $k\in I_{u}$ and $\{f_{p}(Qk+u)\}\in [0,1/p)$ then
\[\left\{\alpha_{j}p^{j-1}(Qk+u)^{j}+\cdots+\alpha_{1}(Qk+u)+\alpha_{0}/p\right\}\in J_{I,u,p}.\]
We will call the above polynomial $g_{p,I,u}(k)$. Let $\mathcal{N}_{I,u}$ be the number of $(p,k)$ with $k\in I_{u}$, $p\in [X,2X]$, $p^{j-1}Q^{j}\alpha_{j}\in \fM_{j}(X,Y_{j},c_{j}\eta)$ and $\{g_{p,I,u}(k)\}\in J_{I,u,p}$ and note that $\mathcal{N}_{I}\le \sum_{u=0}^{Q-1}\mathcal{N}_{I,u}$, because $p^{j-1}\alpha_{j}\in \fM_{j}(X,Y,c_{j}\eta)$ implies that $p^{j-1}Q^{j}\alpha_{j}\in \fM_{j}(X,Y_{j},c_{j}\eta)$. The leading coefficient of $g_{p,I,u}$ is $p^{j-1}Q^{j}\alpha_{j}$.

We first consider the case $j=1$. In this case $g_{p,I,u}(k)=Q\alpha_{1}k+\alpha_{1}u+\alpha_{0}/p$. We must have $\alpha_{1}\in \fM_{1}(X,Y,c_{1}\eta)$ because otherwise $D_{1}(X,Y,c_{1}\eta)=0$. It follows that $Q\alpha_{1}\in \fM_{1}(X,Y,c_{1}\eta)$, so there is some $1\le m\le X^{c_{1}\eta}$, $b\in \bZ$ such that $(b,m)=1$ and $|mQ\alpha_{1}-b|\le N^{c_{1}\eta}Y^{-1}$. Let $\sigma'=\sigma+1$ and $Q'=mQ$. We have $\|Q'\alpha_{i}\|\le X^{\eta}\|Q\alpha_{i}\|\le N^{\sigma'\eta}Y^{-i}$ for $2\le i\le d$ and by the induction hypothesis, if $m\le Q$ then \eqref{lower_degree_coefficients_eq5} follows with $\gamma'=\gamma/2$.

So we may assume that $m>Q$. Note that $\|mQ\alpha_{1}\|\le N^{c_{1}\eta}Y_{1}^{-1}$. By Corollary \ref{weyl_corollary}, we have
\[\#\{k\in I_{u}:\{g_{p,I,u}(k)\}\in J_{I,u,p}\}\ll Y_{1}m^{-1/2}.\]
Now, \eqref{lower_degree_coefficients_eq5} follows because we assumed that $m\ge Q$. This completes the case $j=1$.

Now suppose that $j\ge 2$. By Lemma \ref{major_arcs_alpha_badapprox}, either $D_{1}(X,Y,\bc \eta)\ll_{d,\eta,\bc} NX^{-\eta}$, or there is some integer $1\le m\le X^{3c_{j}\eta+j-2}$ such that $\|m\alpha_{j}\|\le N^{3c_{j}\eta}Y^{-j}$. In the first case, \eqref{lower_degree_coefficients_eq1} follows immediately, so suppose that $m\ge 1$ is minimal with $\|m\alpha_{j}\|\le N^{3c_{j}\eta}Y^{-j}$. Then $\|mQ^{j}\alpha_{j}\|\le N^{3c_{j}\eta}Q^{j}Y^{-j}\le N^{3c_{j}\eta}Y_{j}^{-j}$. By Lemma \ref{q_saving_bound} and the bound $\mu(J_{I,p,u})\ll X^{-\eta}$,
\begin{equation}\label{lower_degree_coefficients_eq4}\mathcal{N}_{I,u}\ll_{d,\eta,\bc,\theta}Y_{j}X^{1-\eta}+XY_{j}m^{-1/(2j^{2})}(\log X)^{-1}.\end{equation}
If $m\ge X^{2j^{2} \eta}$ then this is at most $O(Y_{j}X^{1-\eta})$. Summing over the $Q$ values of $u$ and $O(Y/(QY_{j}))$ values of $I$, we get \eqref{lower_degree_coefficients_eq1}.

So we may assume that $m<X^{2j^{2}\eta}$. In this case, let $\sigma'=2j^{2}+j\sigma$ and $Q'=mQ^{j}\le X^{\sigma'\eta}$. We also have $\|Q'\alpha_{i}\|\le mQ^{j-1}\|Q\alpha_{i}\|\le N^{\sigma'\eta}Y^{-i}$ for $j+1\le i\le d$. If $m\le Q$, then by the induction hypothesis, we have
\[D_{1}(X,Y,\bc \eta)\ll_{d,\bc,\eta} NX^{-\eta}+NQ^{-\gamma}\le NX^{-\eta}+N(Q')^{-\gamma/(j+1)}(\log X)^{-1}.\]
If $m>Q$ then by \eqref{lower_degree_coefficients_eq4}, we have
\[D_{1}(X,Y,\bc \eta)\ll_{f,\bc,\eta} NX^{-\eta}+N(Q')^{1/(2j^{2}(j+1))}(\log X)^{-1}.\]
So if we let $\gamma'=\min\{\gamma/(j+1),1/(2j^{2}(j+1))\}$, then \eqref{lower_degree_coefficients_eq5} follows.
\end{proof}
\end{lemma}

There are two more lemmas that we need before we can bound $\#\mathcal{S}_{1}(N,\theta,\bc\eta)$. From Lemma \ref{lower_degree_coefficients}, we can split $f$ into a part with rational coefficients of small denominator and a part with small coefficients. Lemma \ref{polynomial_small} will allow us to control the part with small coefficients, while Lemma \ref{mod_polynomial_zeros} will allow us to control the part with rational coefficients.

\begin{lemma}\label{polynomial_small}
Let $d\ge 1$ and $h(x)=\sum_{j=0}^{d} \beta_{j}x^{j}$ be a polynomial with $\beta_{j}\in \bR$ and $\beta_{d}\ne 0$. Let $B=\max_{0\le j\le d}|\beta_{j}|$ and for $r\ge 1$, let
\[I(r)=\{x\in [0,1]: |h(x)|\le Br^{-1}\}.\]
Then $I(r)$ is a union of at most $O(d)$ closed intervals and
\[\mu(I(r))\ll_{d} r^{-\frac{2}{d(d+1)}}.\]
The implied constant depends only on $d$.
\begin{proof}
By replacing $h$ with $B^{-1}h$, we may suppose without loss of generality that $B=1$. To see that $I(r)$ is a union of intervals, let $x_{1}<x_{2}\cdots<x_{m}$ be the solutions to $h(x)=\pm r$. For each open interval $(x_{i},x_{i+1})$, either $|h(x)|>r$ for all $x\in (x_{i},x_{i+1})$ or $|h(x)|<r$ for all $x\in (x_{i},x_{i+1})$. Also, since $|h(x)|\rightarrow \infty$ as $|x|\rightarrow \infty$, we have $|h(x)|>r$ when $x>x_{m}$ or $x<x_{1}$.

To prove the upper bound on $\mu(I(r))$, we proceed by induction on $d$. Suppose first that $d=1$. If $|\beta_{1}|< 1/2$, then since $\max_{0\le j\le d}|\beta_{j}|=1$, we must have $|\beta_{0}|=1$ and then $|h(x)|\ge 1/2$ for all $x\in [0,1]$. So if $r\ge 2$ then $I_{r}=\emptyset$ and otherwise $\mu(I_{r})\le 1\ll r^{-1}$. Hence we may suppose that $|\beta_{1}|\ge 1/2$, but then $|h(x)|\le r^{-1}$ implies that $|x-\beta_{0}/\beta_{1}|\le r^{-1}|\beta_{1}|^{-1}$, so $\mu(I_{r})\le 2r^{-1}$.

Now suppose that the lemma holds for $h$ of degree $d-1$. Suppose that $|\beta_{d}|<r^{-(d-1)/(d+1)}$. Let $h^{\ast}(x)=\sum_{j=0}^{d-1}\beta_{j}x^{j}$. Then for all $x\in \bC$,
\[|h(x)|\ge |h^{\ast}(x)|-r^{-\frac{d-1}{d+1}}.\]
So if $|h^{\ast}(x)|\ge 2r^{-(d-1)/(d+1)}$, then $|h(x)|\ge r^{-1}$. We may apply the induction hypothesis to $h^{\ast}$, noting that $\max_{0\le j\le d-1} |\beta_{j}|=1$, to find that
\[\mu(I_{r})\ll_{d} \left(r^{-\frac{d-1}{d+1}}\right)^{\frac{2}{d(d-1)}}=r^{-\frac{2}{d(d+1)}},\]
as required.

Hence we may suppose that $|\beta_{d}|\ge r^{-(d-1)/(d+1)}$. Then there are some $\zeta_{1},\dots,\zeta_{d}\in \bC$ such that
\[h(x)=\beta_{d}\prod_{j=1}^{d}(x-\zeta_{j}).\]
If $|h(x)|<r^{-1}$ then there must be some $1\le j\le d$ such that
\[|x-\zeta_{j}|<\left(|\beta_{d}|^{-1}r^{-1}\right)^{1/d}\le \left(r^{\frac{d-1}{d+1}-1}\right)^{\frac{1}{d}}=r^{-\frac{2}{d(d+1)}}.\]
Hence, $\mu(I_{r})\le 2dr^{-\frac{2}{d(d+1)}}$.
\end{proof}
\end{lemma}

The next lemma is a bound for the number of solutions to a polynomial equation (mod $q$). We will use a similar method to Lemma \ref{mod_kth_roots}.
\begin{lemma}\label{mod_polynomial_zeros}
Let $d\ge 1$. For $q\in \bN$ and a polynomial $g(x)=a_{d}x^{d}+\cdots+a_{0}\in \bZ[x]$ such that $(q,a_{0},\dots,a_{d})=1$, let $\mathcal{C}_{g}(q)$ be the number of solutions $0\le x\le q-1$ to $g(x)\equiv 0$ (mod $q$). Then
\begin{equation}\label{mod_polynomial_zeros_bound}\mathcal{C}_{g}(q)\ll_{d} q^{1-1/(2d^{2})}.\end{equation}
The implied constant depends only on $d$.
\begin{proof}
If $d=1$ then there is either exactly one solution if $(a_{1},q)=1$, or no solutions if $(a_{1},q)>1$ (because $(a_{1},a_{0},q)=1$). For $d\ge2$ we use induction on $d$.

Let $q_{1}=(q,a_{d})$ and $q_{2}=q/q_{1}$. For each solution $0\le x\le q_{1}-1$ to $g(x)\equiv 0$ (mod $q_{1}$), there are at most $q_{2}$ solutions $0\le x\le q-1$ to $g(x)\equiv 0$ (mod $q$). If $g^{\ast}(x)=a_{d-1}x^{d-1}+\cdots+a_{0}$, then $g(x)\equiv g^{\ast}(x)$ (mod $q_{1}$) for all $x$. By the induction hypothesis, we therefore have
\[\mathcal{C}_{g}(q)\le q_{2}\mathcal{C}_{g^{\ast}}(q_{1})\ll_{d} q_{2}q_{1}^{1-1/(2d^{2})}.\]
If $q_{1}\ge q^{1/2}$ then \eqref{mod_polynomial_zeros_bound} follows, so suppose that $q_{1}<q^{1/2}$ from now on. This implies that $q_{2}>q^{1/2}$.

Recall the identity,
\[\sum_{a=0}^{q-1}e(ax/q)=\begin{cases}q,\quad \text{if }q\mid x\\ 0,\quad \text{otherwise}.\end{cases}\]
It follows that
\[\mathcal{C}_{g}(q)=\sum_{x=0}^{q-1} \frac{1}{q} \sum_{y=0}^{q-1} e(yg(x)/q)=\frac{1}{q}\sum_{y=0}^{q-1}\left(\sum_{x=0}^{q-1}e(yg(x)/q)\right).\]
For any $y\in \bZ$, the leading coefficient of $yg(x)/q$ is $ya_{d}/q$, which may be written as $b_{y}/q_{y}$ for some $a_{y}\in \bZ$, where $q_{y}=q_{2}/(y,q_{2})$ and $(a_{y},q_{y})=1$. So by Lemma \ref{weyl_sum_estimate}, for all $\eps>0$ we have
\[\left|\sum_{x=0}^{q-1}e(yg(x)/q)\right|\ll_{\eps} q^{1+\eps}\left(\frac{q_{2}}{(y,q_{2})}\right)^{-\frac{1}{d(d-1)}}.\]
Recalling \eqref{gcd_q_sum}, we have
\[\sum_{y=0}^{q-1}(y,q_{2})^{\frac{1}{d(d-1)}}=\frac{q}{q_{2}}\sum_{y=0}^{q_{2}-1}(y,q_{2})^{\frac{1}{d(d-1)}}\ll_{\eps} q q_{2}^{\eps}.\]
So
\[\frac{1}{q}\sum_{y=0}^{q-1}\sum_{x=0}^{q-1}e(yg(x)/q)\ll_{\eps} q^{\eps}q_{2}^{-\frac{1}{d(d-1)}}\sum_{y=0}^{q-1}(y,q_{2})^{2^{1-d}}\ll_{\eps} q^{1+\eps}q_{2}^{\eps-\frac{1}{d(d-1)}}.\]
Now, recalling that $q_{2}>q^{1/2}$, the result follows.
\end{proof}
\end{lemma}

We are now ready to bound $\#\mathcal{S}_{1}(N,\theta,\bc\eta)$.

\begin{lemma}\label{major_arcs_alpha_goodapprox_med_X}
Suppose that $N\in \bZ$ with $N\ge 1$. Let $\theta\in (0,1/2)$ and $\eta>0$ with $\eta$ sufficiently small in terms of $d,\theta$. Suppose that $j$ is such that $\alpha_{j}$ is irrational and $\alpha_{i}$ is rational for $j+1\le i\le d$. Let $q\in \bN$ be the smallest integer such that $\|q\alpha_{j}\|\le N^{-\theta/2}$. Then there exists $\gamma=\gamma(d,\theta)>0$ such that,
\[\#\mathcal{S}_{1}(N,\theta,\bc\eta)\ll_{f,\eta,\bc,\theta} Nq^{-\gamma}+NR^{-\eta/(2d+1)}.\]
\begin{proof}
By Corollary \ref{major_arcs_alpha_goodapprox_large_X} with $\lambda=1/(2d+1)$, it suffices to prove that
\begin{equation}\label{major_arcs_alpha_goodapprox_med_X_eq1}\#\bigcup_{R\le X\le N^{\frac{1}{2d+1}}}\mathcal{D}_{1}(X,N/X,\bc\eta)\ll_{f,\eta,\bc,\theta} Nq^{-\gamma}+NR^{-\eta},\end{equation}
for some $\gamma>0$.

Let $\sigma>0$ be as in Lemma \ref{lower_degree_coefficients} and consider some $R\le X\le N^{1/(2d+1)}$, $Y=N/X$. Let $Q_{X}\in \bN$ be minimal such that $\|Q_{X}\alpha_{i}\|\le N^{\sigma\eta}Y^{-i}$ for $1\le i\le d$. Then $Y\ge X^{2d}$, so we may apply Lemma \ref{lower_degree_coefficients} to find that $D_{1}(X,Y,\bc\eta)\ll NX^{-\eta}+NQ_{X}^{-\gamma}(\log X)^{-1}$ for some $\gamma>0$ depending only on $d,\bc$. If we let $Q\in \bN$ be minimal such that $\|Q\alpha_{i}\|\le N^{1/2-i}$ for $1\le i\le d$ then we claim that $Q\le Q_{X}$ for all $X\le N^{1/(2d+1)}$ as long as $\eta$ is sufficiently small. This holds because $Y\ge N^{1-1/(2d+1)}$, so $N^{\sigma\eta}Y^{-i}\le N^{\sigma\eta+i/(2d+1)-i}\le N^{1/2-i}$ if $\sigma\eta$ is small enough. It follows that
\begin{equation}\label{major_arcs_alpha_goodapprox_med_X_eq2}D_{1}(X,Y,\bc\eta)\ll NX^{-\eta}+NQ^{-\gamma}(\log X)^{-1},\end{equation}
which suffices when $Q\ge X^{\gamma^{-1}\eta}$. We also have $q\le Q$ by the definition of $q$.

We may therefore suppose from now on that $q\le Q\le X^{\gamma^{-1}\eta}$. Let $a_{1},\dots,a_{d}\in \bZ$ be such that $|Q\alpha_{j}-a_{j}|=\|Q\alpha_{j}\|$ for each $1\le j\le d$. Let
\[E_{1}=\{n\in [N]: Q\mid (a_{d}n^{d}+\cdots+a_{1}n)\}.\]
Note that $E_{1}$ does not depend on $X$, and also that by the minimality of $Q$, we have $(Q,a_{1},\dots,a_{d})=1$.

For each $r\in \bZ$, let $I_{r}$ be the set of $n\in [N]$ such that
\[\frac{a_{d}n^{d}+\cdots+a_{1}n}{Q}+r \le f(n)\le \frac{a_{d}n^{d}+\cdots+a_{1}n}{Q}+r+1.\]
Let $\beta_{j}=\alpha_{j}-a_{j}/Q$ for $1\le j\le d$ and $\beta_{0}=\alpha_{0}$. Recalling that the implied constant in \eqref{major_arcs_alpha_goodapprox_med_X_eq1} may depend on $f$, we may take $N$ to be large enough so that $|\beta_{0}|N^{-1}\le N^{-1/2}$. Let
\[\delta=\max_{0\le j\le d}|N^{j-1}\beta_{j}|\le N^{-1/2}<1.\]
We find that $I_{r}=\emptyset$ unless $|r|\ll \delta N$. We further let $g(n)=\beta_{d}n^{d}+\cdots+\beta_{1}n+\beta_{0}$. Let $E_{2}=\{n\in [N]: |g'(n)|\le \delta Q^{-1}\}$. Applying Lemma \ref{polynomial_small} to the polynomial $g'(Nx)/N$, we see that $E_{2}$ is the intersection of $\bZ$ with a union of at most $O(d)$ intervals of combined length at most $O(NQ^{-\frac{2}{d(d+1)}})$, and it does not depend on $X$ because $g$ does not depend on $X$ (in particular, $\# E_{2}\ll_{d} NQ^{-\frac{2}{d(d+1)}}$). Then $I_{r}\setminus E_{2}$ is the intersection of $\bZ$ with a union of at most $O(d)$ intervals. On each of these intervals, $|g'(x)|> \delta Q^{-1}$, so by an application of the mean value theorem, $|I_{r}\setminus E_{2}|\ll_{d} Q\delta^{-1}$.

Let $D^{\ast}_{1}(X,Y,\bc\eta)$ be the number of $(p,k)$ with $X\le p\le 2X$ prime, $k\in [Y]$, $p^{j-1}\alpha_{j}\in \fM_{j}(X,Y,c_{j}\eta)$ for $1\le j\le d$, $p\mid \floor{f(kp)}$ and $pk\notin E_{1}\cup E_{2}$. Then clearly
\begin{equation}\label{modified_dyadic_decomposition}\#\bigcup_{R\le X\le N^{\frac{1}{2d+1}}}\mathcal{D}_{1}(X,Y,\bc\eta)\le \#E_{1}+\#E_{2}+\sum_{R\le X\le N^{\frac{1}{2d+1}}}D^{\ast}_{1}(X,Y,\bc\eta).\end{equation}

Now, for each $t\in \{0,1,\dots,Q-1\}$, let $u_{t}\in \{0,1,\dots,Q-1\}$ be such that $a_{d}t^{d}+\cdots+a_{1}t\equiv u_{t}$ (mod $Q$). If $n\in I_{r}$ with $n\equiv t$ (mod $Q$) then we either have
\[\lfloor f(n)\rfloor = \frac{a_{d}n^{d}+\cdots+a_{1}n-u_{t}}{Q}+r\]
or
\[\lfloor f(n)\rfloor = \frac{a_{d}n^{d}+\cdots+a_{1}n-u_{t}}{Q}+r+1.\]
Suppose that $p\mid n$ and $p\mid \lfloor f(n)\rfloor$. Then $p\mid (a_{d}n^{d}+\cdots+a_{1}n)$, so either $p\mid u_{t}+Qr$ or $p\mid u_{t}+Q(r+1)$. If $n\notin E_{1}$ then $Q\nmid (a_{d}n^{d}+\cdots+a_{1}n)$, so $u_{t}\ne 0$ and $u_{t}+Qr$ and $u_{t}+Q(r+1)$ are both nonzero. In particular, the interval $[Qr, Q(r+2)]$ contains a nonzero multiple of $p$. Recalling that $I_{r}=\emptyset$ unless $|r|\ll \delta N$, and that $Q\le X^{\gamma^{-1}\eta}\le p^{\gamma^{-1}\eta}$, we find that there are at most $O(\delta N Q p^{-1})\le O(\delta N p^{\gamma^{-1}\eta-1})$ values of $r$ such that $[Qr,Q(r+2)]$ contains a nonzero multiple of $p$. For each of these $r$, recall that $I_{r}\setminus E_{2}$ consists of at most $O(d)$ intervals and $|I_{r}\setminus E_{2}|\ll_{d} Q\delta^{-1}$ so the number of $n\in I_{r}\setminus E_{2}$ which are divisible by $p$ is at most $O_{d}(1+\delta^{-1}Qp^{-1})$.

So the number of $n\in [N]\setminus (E_{1}\cup E_{2})$ such that $p\mid n$ and $p\mid \lfloor f(n)\rfloor$ is at most $O_{d}(\delta Np^{\gamma^{-1}\eta-1}+Np^{2\gamma^{-1}\eta-2})$. Recalling that $\delta\le N^{-1/2}$, we see that this is $O_{d}(Np^{2\gamma^{-1}\eta-2}+N^{1/2}p^{\gamma^{-1}\eta-1})$.

By choosing $\eta$ sufficiently small in terms of $d,\bc$, we have $\gamma^{-1}\eta<1/4$ and $X\le N^{1/(2d+1)}\le N^{1/5}$. It follows that
\[D^{\ast}_{1}(X,Y,\bc\eta)\ll_{d} NX^{2\gamma^{-1}\eta-1}+N^{1/2}X^{\sigma\eta}\ll NX^{-1/2}.\]

So by \eqref{modified_dyadic_decomposition}, we have
\begin{align*}\#\bigcup_{R\le X\le N^{\frac{1}{2d+1}}}\mathcal{D}_{1}(X,Y,\bc\eta)&\le \#E_{1}+\#E_{2}+\sum_{R\le X\le N^{\frac{1}{2d+1}}}NX^{-1/2}\\ &\ll \#E_{1}+\#E_{2}+R^{-1/2}.\end{align*}
Recall that $\#E_{2}\ll NQ^{-\frac{2}{d(d+1)}}\le Nq^{-\frac{2}{d(d+1)}}$. By Lemma \ref{mod_polynomial_zeros}, we have $\#E_{1}\ll_{d} NQ^{-1/(2d^{2})}\le Nq^{-1/(2d^{2})}$.

The result follows.
\end{proof}
\end{lemma}

\section{Proof of Theorem \ref{density_zero}}
Putting together the bounds on $\#\mathcal{S}_{1}$, $\#\mathcal{S}_{2}$ and $\#\mathcal{S}_{3}$ from the previous sections, we have the following.
\begin{lemma}\label{calD_bound}
Let $d\in \bZ$ with $d\ge2$ and $f(x)=\alpha_{d} x^{d}+\alpha_{d-1} x^{d-1}+\cdots+\alpha_{0}$, with $\alpha_{i}\in \bR$ for $0\le i\le d$. Suppose that $1\le j\le d$ is such that $\alpha_{j}$ is irrational and $\alpha_{i}$ is rational for $j+1\le i\le d$. Let $\theta\in (0,1/2)$ and $\eta>0$ sufficiently small in terms of $\theta,d$ and let $q\in \bN$ be minimal such that $\|q\alpha_{j}\|\le N^{-\theta/2}$. Then there is some $C(f,\theta,\eta)>0$ such that for all $R\ge 3$,
\[\#\mathcal{S}_{f,R}(N)\ll \theta N+\frac{N}{\log N}+C(f,\theta,\eta)N\left(q^{-\gamma}+R^{-\eta/(2d^{2})}\right).\]
\begin{proof}
Lemma \ref{minor_arc_bound} tells us that there is some $\bc=(c_{1},\dots,c_{d})\in \bR_{\ge 1}^{d}$ (depending on $d,\theta$) such that
\[\#\mathcal{S}_{2}(N,\theta,\bc \eta)\ll_{d,\eta,\theta} NR^{-\eta/(2d^{2})}.\]
Then applying Lemma \ref{major_arcs_alpha_goodapprox_med_X} with this choice of $\bc$, we have some $\gamma=\gamma(d,\theta)>0$ such that
\[\#\mathcal{S}_{1}(N,\theta,\bc\eta)\ll_{f,\eta,\theta} Nq^{-\gamma}+NR^{-\eta/(2d+1)}.\]
Then recalling \eqref{calS_arc_decomposition}, we have
\[\#\mathcal{S}_{f,R}(N)\ll  \theta N+\frac{N}{\log N}+ \#\mathcal{S}_{1}(N,\theta,\boldsymbol{\eta})+\#\mathcal{S}_{2}(N,\theta,\boldsymbol{\eta}),\]
so there is some constant $C(f,\theta,\eta)>0$ such that
\[\#\mathcal{S}_{f,R}(N)\ll \theta N+\frac{N}{\log N}+C(f,\theta,\eta)N\left(q^{-\gamma}+R^{-\eta/(2d^{2})}\right),\]
as required.
\end{proof}
\end{lemma}

We can now finally deduce Theorem \ref{density_zero}.

\begin{proof}[Proof of Theorem \ref{density_zero}]
By definition,
\[\overline{\dens}\{n\in \bN: \exists p\ge R \: \text{prime},\: p\mid n\:\text{and}\: p\mid \lfloor f(n)\rfloor\}=\lim_{N\rightarrow\infty} \frac{1}{N}\#\mathcal{S}_{f,R}(N).\]
Let $\theta>0$. There is some $1\le j\le d$ such that $\alpha_{j}$ is irrational, and we may suppose that $j$ is maximal with this property so that $\alpha_{i}$ is rational for all $i\ge j+1$. Let $q=q(N,\theta)\in \bN$ be minimal such that $\|q\alpha_{j}\|\le N^{-\theta/2}$. Note that for any fixed $\theta>0$, we must have $q(N,\theta)\rightarrow \infty$ as $N\rightarrow \infty$. Hence, by Lemma \ref{calD_bound} we have some $\eta>0$ and $C(f,\theta,\eta)>0$ such that
\[\limsup_{N\rightarrow\infty} \frac{1}{N}\#\mathcal{S}_{f,R}(N)\ll \theta+C(d,\theta)R^{-\eta/(2d^{2})}.\]
So,
\[\limsup_{R\rightarrow\infty}\limsup_{N\rightarrow\infty} \frac{1}{N}\#\mathcal{S}_{f,R}(N)\ll \theta.\]
But since $\theta>0$ is arbitrary, the limit is in fact $0$.
\end{proof}

\section{Acknowledgements}
The author was supported by the National Research, Development and Innovation Office NKFIH (Excellence program, Grant Nr. 153829).

\end{document}